\documentclass[11pt]{article}
\usepackage{etoolbox}
\usepackage{amsmath,amsthm,amsfonts,amssymb}
\usepackage[numbers,sort&compress]{natbib}
\usepackage{color,thmtools,thm-restate}
\usepackage[colorlinks=true, citecolor=black, linkcolor=black,
urlcolor=blue,anchorcolor=blue,citecolor=black, filecolor=blue,
menucolor=blue, pagecolor=blue, pdfcreator=Me, pdfproducer=Me]{hyperref}
\usepackage{microtype}
\usepackage[margin=30mm]{geometry}
\usepackage{enumerate,calc,graphicx,float}
\usepackage[T1]{fontenc}
\usepackage{soul}

\usepackage[noabbrev,capitalise]{cleveref}
\crefname{lem}{Lemma}{Lemmas}
\crefname{thm}{Theorem}{Theorems}
\crefname{prop}{Proposition}{Propositions}

\theoremstyle{plain}
\newtheorem{thm}{Theorem}
\newtheorem{lem}[thm]{Lemma}
\newtheorem{cor}[thm]{Corollary}
\newtheorem{prop}[thm]{Proposition}
\newtheorem{conj}[thm]{Conjecture}

\renewcommand{\baselinestretch}{1.19}
\usepackage[hang]{footmisc} 

\renewcommand{\thefootnote}{\fnsymbol{footnote}}	
\allowdisplaybreaks
\newcommand{\arXiv}[1]{arXiv:\,\href{http://arxiv.org/abs/#1}{#1}}

\newcommand{\MSN}[2]{MR:\,\href{http://www.ams.org/mathscinet-getitem?mr=MR#1}{#1}}
\newcommand{\doi}[1]{doi:\,\href{http://dx.doi.org/#1}{#1}}

\newcommand{\website}[1]{\href{http://#1}{#1}}

\newcommand\DateFootnote{
  \begingroup
  \renewcommand\thefootnote{}
  \footnote{21st April 2017, revised 23rd February 2018\\
MSC:  05C15, 05C83}  
  \setcounter{footnote}{0}
  \vspace*{-3ex}
  \endgroup}

\makeatletter
\renewcommand\section{\@startsection {section}{1}{\z@}%
  {-3ex \@plus -1ex \@minus -.2ex}%
  {2ex \@plus.2ex}%
  {\normalfont\large\bfseries}}
\renewcommand\subsection{\@startsection{subsection}{2}{\z@}%
  {-2.5ex\@plus -1ex \@minus -.2ex}%
  {1.5ex \@plus .2ex}%
  {\normalfont\normalsize\bfseries}}
\renewcommand\subsubsection{\@startsection{subsubsection}{3}{\z@}%
  {-2ex\@plus -1ex \@minus -.2ex}%
  {1ex \@plus .2ex}%
  {\normalfont\normalsize\bfseries}}
\renewcommand\paragraph{\@startsection{paragraph}{4}{\z@}%
  {1.5ex \@plus.5ex \@minus.2ex}%
  {-1em}%
  {\normalfont\normalsize\bfseries}}
\renewcommand\subparagraph{\@startsection{subparagraph}{5}{\parindent}%
  {1.5ex \@plus.5ex \@minus .2ex}%
  {-1em}%
  {\normalfont\normalsize\bfseries}}
\makeatother

\renewcommand{\thefootnote}{\fnsymbol{footnote}}

\newcommand{\bceil}[1]{\bigl\lceil{#1}\bigr\rceil}

\newcommand{\half}{\ensuremath{\protect\tfrac{1}{2}}}
\renewcommand{\geq}{\geqslant}
\renewcommand{\leq}{\leqslant}
\renewcommand{\preceq}{\preccurlyeq}
\newcommand{\Oh}[1]{\ensuremath{\protect\mathcal{O}(#1)}}
\newcommand{\preproof}{\vspace*{-3ex}}

\DeclareMathOperator{\scol}{scol}
\DeclareMathOperator{\wcol}{wcol}
\newcommand{\wreach}{W}
\newcommand{\sreach}{S}
\DeclareMathOperator{\dist}{dist}

\title{\Large\bfseries Improper Colourings inspired by Hadwiger's
  Conjecture\,\footnotemark[1]}
\author{Jan van den Heuvel\,\footnotemark[2] \ \ and \ \
  David~R.~Wood\,\footnotemark[3]}
\date{}

\begin{document}
\maketitle

\DateFootnote

\footnotetext[1]{Research for this paper was done during a visit of the
  first author to Monash University. JvdH would like to thank the School of
  Mathematical Sciences at Monash University for hospitality and support.}

\footnotetext[2]{Department of Mathematics, London School of Economics and
  Political Science, United Kingdom (\texttt{j.van-den-heuvel@lse.ac.uk}).}

\footnotetext[3]{School of Mathematical Sciences, Monash University,
  Melbourne, Australia (\texttt{david.wood@monash.edu}). Supported by the
  Australian Research Council.}

\begin{abstract}
  \noindent
  Hadwiger's Conjecture asserts that every $K_t$-minor-free graph has a
  proper $(t-1)$-colouring. We relax the conclusion in Hadwiger's
  Conjecture via improper colourings. We prove that every $K_t$-minor-free
  graph is $(2t-2)$-colourable with monochromatic components of order at
  most $\bceil{\half(t-2)}$. This result has no more colours and much
  smaller monochromatic components than all previous results in this
  direction. We then prove that every $K_t$-minor-free graph is
  $(t-1)$-colourable with monochromatic degree at most $t-2$. This is the
  best known degree bound for such a result. Both these theorems are based
  on a decomposition method of independent interest. We give analogous
  results for \mbox{$K_{s,t}$-minor}-free graphs, which lead to improved
  bounds on generalised colouring numbers for these classes. Finally, we
  prove that graphs containing no $K_t$-immersion are $2$-colourable with
  bounded monochromatic degree.
\end{abstract}

\renewcommand{\thefootnote}{\arabic{footnote}}

\section{Introduction}
\label{Introduction}

Hadwiger's Conjecture \citep{H43} asserts that every $K_t$-minor-free graph
has a proper $(t-1)$-colouring. For $t\leq 3$ the conjecture is easy.
\citet{H43} and \citet{Dirac52} independently proved the conjecture for
$t=4$; while Wagner's result \citep{Wagner37} means that the case $t=5$ is
equivalent to the Four Colour Theorem. Finally, \citet{RST-Comb93} proved
Hadwiger's Conjecture for $t=6$. The conjecture remains open for $t\geq7$.
Hadwiger's Conjecture is widely considered to be one of the most important
open problems in graph theory. The best upper bound on the chromatic number
of $K_t$-minor-free graphs is $\Oh{t\sqrt{\log t}}$ independently due to
Kostochka~\citep{Kostochka82,Kostochka84} and
\citet{Thomason84,Thomason01}. See the recent survey by \citet{SeymourHC}
for more on Hadwiger's Conjecture.

One possible way to approach Hadwiger's Conjecture is to allow improper
colourings. In a vertex-coloured graph, a \emph{monochromatic component} is
a connected component of the subgraph induced by all the vertices of one
colour. A graph $G$ is \emph{$k$-colourable with clustering $c$} if each
vertex can be assigned one of $k$ colours such that each monochromatic
subgraph has at most~$c$ vertices\footnote{This type of colouring is
  sometimes called ``fragmented'' in the literature, but we feel that
  increased fragmentation suggests smaller monochromatic components, hence
  we use the term ``clustered''.}. \citet{KMRV97} introduced this type of
colouring, and now many results are known. The \emph{clustered chromatic
  number} of a graph class $\mathcal{G}$ is the minimum integer $k$ for
which there exists an integer $c$ such that every graph in $\mathcal{G}$ is
$k$-colourable with clustering $c$.

\citet{KawaMohar-JCTB07} were the first to prove an $\Oh{t}$ upper bound on
the clustered chromatic number of $K_t$-minor-free graphs. In particular,
they proved that every $K_t$-minor-free graph is
$\bceil{\frac{31}{2}t}$-colourable with clustering $f(t)$, for some
function $f$. The number of colours in this result was improved to
$\bceil{\half(7t-3)}$ by \citet{Wood10}\,\footnote{This result depends on a
  result announced in 2008 which is not yet written.}, to $4t-4$ by
\citet{EKKOS}, and to $3t-3$ by \citet{LiuOum}. See \citep{Kawa08,KO16} for
analogous results for graphs excluding odd minors. For all of these
results, the function $f(t)$ is very large, often depending on constants
from the Graph Minor Structure Theorem \citep{RS99}.

Our first contribution is to prove an analogous theorem with the best known
number of colours, and also with small clustering. The proof is simple, and
does not depend on any deep theory.

\begin{thm}
  \label{MonoComponents}\mbox{}\\*
  For $t\geq4$, every $K_t$-minor-free graph is $(2t-2)$-colourable with
  clustering $\bceil{\half(t-2)}$.
\end{thm}

\cref{MonoComponents} implies that the clustered chromatic number of
$K_t$-minor-free graphs is at most $2t-2$. A construction of Edwards et
al.~\cite{EKKOS} mentioned below implies that the clustered chromatic
number of $K_t$-minor-free graphs is at least $t-1$.

A second way to relax the conclusion in Hadwiger's Conjecture is to bound
the maximum degree of monochromatic components. A graph $G$ is
\emph{$k$-colourable with defect $d$} if each vertex can be assigned one of
$k$ colours such that each vertex is adjacent to at most $d$ vertices of
the same colour; that is, each monochromatic subgraph has maximum degree at
most $d$. \citet{CCW86} introduced the notion of defective graph colouring,
and now many results for various graph classes are known. A graph class
$\mathcal{G}$ is \emph{defectively $k$-colourable} if there exists an
integer~$d$ such that every graph in $\mathcal{G}$ is $k$-colourable with
defect $d$. The \emph{defective chromatic number} of~$\mathcal{G}$ is the
minimum integer $k$ such that $\mathcal{G}$ is defectively $k$-colourable
\citep{CGJ97}. \citet{EKKOS} proved that every $K_t$-minor-free graph is
$(t-1)$-colourable with defect $\Oh{t^2\log t}$. Moreover, it is shown in
\citep{EKKOS} that the number of colours, $t-1$, is best possible in the
following strong sense: for every integer $d$, there is a $K_t$-minor-free
graph that is not $(t-2)$-colourable with defect $d$. Thus the defective
chromatic number of $K_t$-minor-free graphs equals $t-1$. (This also shows
that the clustered chromatic number of $K_t$-minor-free graphs is at least
$t-1$.)

Our second contribution is an improved upper bound on the defect in the
result of \citet{EKKOS}.
 
\begin{thm}
  \label{MaxDegree}\mbox{}\\*
  For $t\geq4$, every $K_t$-minor-free graph is $(t-1)$-colourable with
  defect $t-2$.
\end{thm}

\citet{EKKOS} wisely noted that their theorem mentioned earlier should not
be considered evidence for the truth of Hadwiger's Conjecture, since their
method also proves that every \mbox{$K_t$-topological}-minor-free graph is
$(t-1)$-colourable with defect $\Oh{t^4}$. It is not true that every
$K_t$-topological-minor-free graph is properly $(t-1)$-colourable. This
\emph{last statement} is Haj\'os' Conjecture, which is now known to be
false \citep{Catlin79,Thomassen05}. On the other hand, our proof of
\cref{MaxDegree} does not work for graphs excluding a topological minor.

\cref{MonoComponents,MaxDegree} are corollaries of the following
decomposition result of independent interest. A sequence $H_1,\dots,H_\ell$
is a \emph{connected partition} of a graph $G$ if each $H_i$ is a non-empty
connected induced subgraph of $G$, the subgraphs $H_1,\dots,H_\ell$ are
pairwise disjoint, and $V(G)=V(H_1)\cup\dots\cup V(H_\ell)$. Two disjoint
subgraphs $H$ and $H'$ of a graph $G$ are \emph{adjacent} if there is an
edge in $G$ with one endpoint in $H$ and one endpoint in $H'$. For a
positive integers $n,m$, we use $[n]$ to denote the set $\{1,\ldots,n\}$
and $[n,m]$ to denote the set $\{n,\ldots,m\}$.

\begin{thm}
  \label{BasicDecomposition}\mbox{}\\*
  For $t\geq4$, every $K_t$-minor-free graph $G$ has a connected partition
  $H_1,\dots,H_\ell$ such that for $i\in[\ell]$:
  \vspace*{-2.25ex}
  \begin{enumerate}[(1)]
    \setlength\itemsep{0ex}
  \item $H_i$ is adjacent to at most $t-2$ of the subgraphs
    $H_1,\dots,H_{i-1}$;
  \item $H_i$ has maximum degree at most $t-2$; and
  \item $H_i$ is $2$-colourable with clustering $\bceil{\half(t-2)}$.
  \end{enumerate}
\end{thm}

We actually prove a decomposition theorem with several further properties;
see \cref{Decomposition}. It is easy to derive
\cref{MaxDegree,MonoComponents} from \cref{BasicDecomposition}. Colour the
subgraphs $H_1,\dots,H_\ell$ greedily in this order, such that adjacent
subgraphs receive distinct colours. By property (1), $t-1$ colours suffice.
\cref{MaxDegree} follows from property (2) by colouring each vertex in
$H_i$ by the colour assigned to $H_i$. \cref{MonoComponents} follows from
property (3) by taking the product of the $(t-1)$-colouring of
$H_1,\dots,H_\ell$ with the given $2$-colouring of each subgraph $H_i$.

\cref{BasicDecomposition} is an extension of a result by \citet{HOQRS} in
which properties (2) and~(3) are replaced by ``$H_i$ has a Breadth-First
Search (BFS) spanning tree with at most $t-3$ leaves''. \citet{HOQRS} were
motivated by connections to generalised colouring numbers. Note that the
result in \citep{HOQRS} implies that~$H_i$ has at most $t-3$ vertices in
each BFS layer. It follows that the maximum degree of~$H_i$ is at most
$3t-10$. Alternately colouring the BFS layers shows that~$H_i$ is
\mbox{$2$-colourable} with clustering $t-3$. Constructing $H_i$ more
carefully, and choosing the $2$-colouring more carefully, leads to the
improved bounds in \cref{BasicDecomposition}, which we prove in
\cref{ProofDecomposition}.

Our main decomposition theorem, \cref{Decomposition}, also has the
following corollary, which might be of independent interest.

\begin{thm}\label{TreeWidth}\mbox{}\\*
  For $t\geq4$, every $K_t$-minor-free graph $G$ has a connected partition
  $H_1,\dots,H_\ell$ such that:
  \vspace*{-2.25ex}
  \begin{enumerate}[(1)]
    \setlength\itemsep{0ex}
  \item the quotient graph~$Q$ obtained by contracting each $H_i$ to a
    single vertex is chordal with clique size at most $t-1$ (and hence has
    treewidth at most $t-2$); and
  \item each part $H_i$ has bandwidth (and hence pathwidth and treewidth)
    at most $t-3$.
  \end{enumerate}
\end{thm}

Hadwiger's Conjecture implies that for every graph $H$ with $t$ vertices,
the maximum chromatic number of $H$-minor-free graphs equals $t-1$ (since
$K_{t-1}$ is $H$-minor-free). However, for clustered and defective
colourings, fewer colours often suffice. For example, it follows from the
main result by \citet{OOW16} that for every fixed non-complete graph $H$
on~$t$ vertices, every $H$-minor-free graph is $(t-2)$-colourable with
bounded defect, which is one fewer colour than in the complete graph case.
More interestingly,, \citet{Archdeacon1987} proved that graphs embeddable
in a fixed surface are defectively $3$-colourable (see also
\citep{Woodall11,CE16,CGJ97,CCW86}); while \citet{DN17} proved that such
graphs are $4$-colourable with bounded clustering. \citet{OOW16}
conjectured that for every connected graph $H$, the defective chromatic
number of $H$-minor-free graphs equals the treedepth of $H$ minus $1$. They
proved this conjecture for $K_{s,t}$-minor-free graphs, by showing that
they are defectively \mbox{$s$-colourable}. Note that $K_{3,t}$-minor-free
graphs are of particular interest since they include and generalise graphs
embeddable in fixed surfaces. In the case $s\leq3$ we prove decomposition
results analogous to \cref{BasicDecomposition} that imply this result of
\citep{OOW16} with much improved bounds on the defect. This direction is
explored in \cref{K2tK3t}.

In the same way as \citet{HOQRS} applied their decomposition result for
\mbox{$K_t$-minor}-free graphs to the setting of generalised colouring
numbers, we apply our decomposition results for $K_{s,t}$-minor-free graphs
and $K^*_{s,t}$-minor-free graphs (where $K^*_{s,t}$ is the complete join
of~$K_s$ and $\overline{K_t}$) to conclude new bounds on generalised
colouring numbers. Our results when specialised for graphs of given genus
are almost as strong as the best known bounds. We then show how the concept
of layered treewidth also leads to good bounds on generalised colouring
numbers. The advantage of this approach is that it immediately applies to
several non-minor-closed graph classes of interest. These results on
generalised colouring numbers are presented in \cref{GeneralisedColouring}.

The final section, \cref{Immersions}, returns to the topic of defective
graph colouring, but instead of excluding a $K_t$ minor we exclude a $K_t$
immersion. The analogue of Hadwiger's Conjecture, that $K_t$-immersion-free
graphs are properly $(t-1)$-colourable \citep{LM89,AL03}, is open. For
defective colouring, we show that only 2 colours suffice.

Before continuing, we mention an important connection between clustered and
defective colourings, implicitly observed in \citep{EKKOS}. We include the
proof for completeness.

\begin{lem}[\citet{EKKOS}]
  \label{Defective2Fragmented}\mbox{}\\*
  For every minor-closed graph class $\mathcal{G}$, the clustered chromatic
  number of $\,\mathcal{G}$ is at most three times the defective chromatic
  number of $\,\mathcal{G}$.
\end{lem}

\preproof\begin{proof}
  \citet{LiuOum} proved that for every minor-closed graph class
  $\mathcal{G}$ and integer $d$, there is an integer $c=c(\mathcal{G},d)$
  such that every graph in $\mathcal{G}$ with maximum degree $d$ is
  \mbox{$3$-colourable} with clustering $c$. (\citet{EJ14} previously
  proved an analogous result for graphs on surfaces.) Let $k$ be the
  defective chromatic number of $\mathcal{G}$. Thus for some integer~$d$,
  every graph~$G$ in $\mathcal{G}$ is $k$-colourable with defect $d$. Apply
  the result of Liu and Oum~\citep{LiuOum} to each monochromatic component
  of $G$, which has maximum degree at most $d$. Then $G$ is
  \mbox{$3k$-colourable} with clustering $c$, and hence the clustered
  chromatic number of $\mathcal{G}$ is at most~$3k$.
\end{proof}

\section{Preliminaries}

\subsection{Notation and definitions}

This subsection briefly states standard graph theoretic definitions
probably familiar to most readers.

A graph $H$ is a \emph{minor} of a graph $G$ if a graph isomorphic to $H$
can be obtained from a subgraph of $G$ by contracting edges. Equivalently,
and often easier to use intuitively: a graph~$H$ with vertices
$v_1,\ldots,v_n$ is a minor of $G$ if there exist pairwise disjoint
connected subgraphs $H_1,\ldots,H_n$ of~$G$ such that for every edge
$v_iv_j$ in $H$, $H_i$ and $H_j$ are adjacent in~$G$. We call~$H_i$ the
\emph{branch set corresponding to~$v_i$}. A class of graphs~$\mathcal{G}$
is \emph{minor-closed} if for every graph $G\in\mathcal{G}$, every minor of
$G$ is also in $\mathcal{G}$. A graph $H$ is a \emph{topological minor} of
a graph $G$ if a graph isomorphic to a subdivision of $H$ is a subgraph of
$G$.

The \emph{Euler genus} of an orientable surface with $h$ handles is $2h$.
The \emph{Euler genus} of a non-orientable surface with $c$ cross-caps is
$c$. The \emph{Euler genus} of a graph $G$ is the minimum Euler genus of a
surface in which $G$ embeds (with no crossing edges).

A \emph{tree decomposition} of a graph $G$ is given by a tree $T$ whose
nodes index a collection\linebreak $(\,T_x\subseteq V(G)\mid x\in V(T)\,)$
of sets of vertices in $G$ called \emph{bags}, such that (1) for every edge
$vw$ of $G$, some bag $T_x$ contains both $v$ and $w$, and (2) for every
vertex $v$ of $G$, the set $\{\,x\in V(T)\mid v\in T_x\,\}$ induces a
non-empty (connected) subtree of $T$. The \emph{width} of a tree
decomposition $T$ is $\max\{\,|T_x\bigm|x\in V(T)\,\}-1$, and the
\emph{treewidth} of a graph $G$ is the minimum width of the tree
decompositions of $G$. Note that the treewidth of $G$ equals the minimum
integer $k$ such that $G$ is a subgraph of a chordal graph with clique
number $k+1$.

A \emph{path decomposition} is a tree decomposition in which the underlying
tree is a path. The \emph{pathwidth} of a graph $G$ is the minimum width of
a path decomposition of $G$.

For a graph $G$ and $A,B\subseteq V(G)$, an \emph{$AB$-separator} is a set
$S\subseteq V(G)$ such that every \mbox{$AB$-path} in $G$ contains a vertex
from $S$. (Note that we allow $A$ and $B$ to intersect and that all
vertices in $A\cap B$ must be included in any $AB$-separator.) A pair
$(G_1,G_2)$ is a \emph{$k$-separation} of a graph~$G$ if $G_1$ and $G_2$
are induced subgraphs of $G$ such that $G=G_1\cup G_2$,
$G_1\not\subseteq G_2$ and $G_2\not\subseteq G_1$, and $|V(G_1)\cap
V(G_2)|=k$.

\subsection{Connected Induced Subgraphs}
\label{ConnectedInducedSubgraphs}

This subsection contains some elementary results about connected induced
subgraphs containing a given set of vertices. We look in detail at
so-called \emph{Lexicographic Breadth-First Search} (LexBFS) trees, since
these form a key tool in our methods.

A \emph{layering} of a graph $G$ is a partition $(V_0,V_1,\dots,V_\ell)$ of
$V(G)$ such that for every edge $vw\in E(G)$, if $v\in V_i$ and $w\in V_j$,
then $|i-j|\leq1$. Each set $V_i$ is called a \emph{layer}.

Let $r$ be a vertex in a connected graph $G$. Let
$\ell=\max\{\,\dist_G(r,v)\mid v\in V(G)\,\}$, and for
$i\in\{0,\ldots,\ell\}$ define $V_i=\{\,v\in V(G)\mid\dist_G(r,v)=i\,\}$.
Then $V_0,V_1,\dots,V_\ell$ is a layering of~$G$, called the \emph{BFS
  layering} of $G$ starting from the \emph{root} $r$; each~$V_i$ is called
a \emph{BFS layer} of $G$. A spanning tree $T$ of $G$ rooted at $r$ is a
\emph{BFS spanning tree} if $\dist_G(v,r)=\dist_T(v,r)$ for every vertex
$v$ in $G$. A \emph{BFS subtree} is a subtree of a BFS spanning tree that
includes the root. Let~$S$ be a BFS subtree rooted at $r$ and consider a
vertex $v\in V_i\cap V(S)$ for some $i\geq1$. Let $P_v$ be the $vr$-path in
$S$. Then $P_v$ has exactly one vertex in each of $V_0,\dots,V_i$. The
\emph{parent} of $v$ is the neighbour of $v$ (in $S$) in~$V_{i-1}$. Every
vertex $x$ in $G$ is adjacent to at most three vertices in $P_v$ (since if
$x\in V_j$, then $N_G(x)\subseteq V_{j-1}\cup V_j\cup V_{j+1}$). A
\emph{leaf} in a rooted tree is a non-root vertex of degree 1. If $S$
has~$p$ leaves, then every vertex in $G$ is adjacent to at most $3p$
vertices in $S$. This observation can be improved for a special type of BFS
(sub)trees.

For our purposes, a BFS spanning tree $T$ of $G$ is a \emph{LexBFS spanning
  tree} if each BFS layer~$V_i$ can be linearly ordered such that:
\vspace*{-2.5ex}
\begin{enumerate}[(a)]
  \setlength\itemsep{0ex}
\item each vertex $v\in V_i$ with parent $w\in V_{i-1}$ in $T$ has no
  neighbour in $G$ that comes before $w$ in the ordering of $V_{i-1}$
  (called the \emph{priority rule}); and
\item for every edge $vw$ in $T$ with $v\in V_i$ and $w\in V_{i-1}$, there
  is no edge $xy$ in $T$ with $x$ before~$v$ in the ordering of $V_i$ and
  $y$ after $w$ in the ordering of $V_{i-1}$ (called the \emph{non-crossing
    rule}).
\end{enumerate}
\vspace*{-2ex}
It is easily seen that every connected graph has a LexBFS spanning tree
rooted at any given vertex. A \emph{LexBFS subtree} is a subtree of a
LexBFS spanning tree that includes the root.

Throughout this paper we follow the convention that the root of a rooted
tree (such as a BFS or LexBFS (sub)tree) is never a leaf.

\clearpage
\begin{lem}
  \label{LexBFSDegree}\mbox{}\\*
  For $k\geq1$, if $S$ is a LexBFS subtree of a connected graph $G$ and $S$
  has $k$ leaves, then every vertex in~$G$ has at most $2k$ neighbours in
  $V(S)$.
\end{lem}

\preproof\begin{proof}
  Let $T$ be a LexBFS spanning tree of $G$, such that $S$ is a subtree of
  $T$. Let $V_0,\dots,V_\ell$ be the BFS layers of $T$. Let $v$ be a vertex
  in $V_i$ (which may or may not be in $S$). If $v$ is on some leaf-root
  path $P$ of $S$, then $|N_G(v)\cap V(P)|\leq2$. Now consider a leaf-root
  path $P$ in $S$ not containing $v$. Suppose on the contrary that there
  are distinct vertices $x,y,z\in N_G(v)\cap V(P)$, none of which are on a
  leaf-root path of $T$ containing $v$. Without loss of generality,
  $x\in V_{i-1}$, $y\in V_i$ and $z\in V_{i+1}$. Let $w$ be the parent of
  $v$ in $T$. So $w\in V_{i-1}$, but $w\ne x$ (since $x$ is not on a
  leaf-root path of $T$ containing $v$). By the priority rule, $w$ comes
  before $x$ in the ordering of~$V_{i-1}$. By the non-crossing rule, $v$
  comes before $y$ in the ordering of $V_i$, which contradicts the priority
  rule for $z$. Thus $|N_G(v)\cap V(P)|\leq2$. Since there are $k$
  leaf-root paths in~$S$, in total this gives $|N_G(v)\cap V(S)|\leq2k$.
\end{proof}

A graph $G$ has \emph{bandwidth} at most $k$ if there is a vertex ordering
$v_1,\dots,v_n$ of $V(G)$, such that $|i-j|\leq k$ for each edge $v_iv_j$
of $G$.

\begin{lem}
  \label{LexBFSBandwidth}\mbox{}\\*
  Every connected graph $G$ that has a LexBFS spanning tree $T$ with $k$
  leaves has bandwidth, pathwidth and treewidth at most $k$.
\end{lem}

\preproof\begin{proof}
  Say $T$ is rooted at $r$. Let $V_0,\dots,V_\ell$ be the BFS layers of
  $T$. Each $V_i$ is linearly ordered by LexBFS. We claim that the
  vertex-ordering of $V(G)$ produced by using the orderings of
  $V_0,\dots,V_\ell$ in that order has bandwidth at most $k$. Consider an
  edge $vw$ where $v\in V_i$ and $w\in V_i$. Since $T$ has at most $k$
  leaves, $|V_i|\leq k$ and at most $k-2$ vertices are between~$v$ and $w$
  in $V_0,\dots,V_\ell$. Now consider an edge $vw$ where $v\in V_i$ and
  $w\in V_{i+1}$. Let $X$ be the set of vertices that come after $v$ in
  $V_i$ or come before $w$ in $V_{i+1}$. Then $X$ is the set of vertices
  between~$v$ and $w$ in the ordering of $V(G)$. Let $p$ be the parent of
  $w$ in $T$. By the priority rule, $p\notin X$. By the non-crossing rule,
  no vertex in $X\cup\{v\}$ is a descendent of another vertex in
  $X\cup\{v\}$. Hence, the number of leaves in $T$ is at least
  $|X|+1$, implying $|X|\leq k-1$. Therefore~$G$ has bandwidth at most $k$.

  It is well-known and easy to prove that the pathwidth of a graph is at
  most its bandwidth (and hence so is the treewidth). Take the vertex
  ordering $v_1,\dots,v_n$ of $V(H)$ that shows $H$ has bandwidth at most
  $k$. For $i\in[n-k]$, let $T_i=\{v_i,\ldots,v_{i+k}\}$. Then
  $T_1,T_2,\dots,T_{n-k}$ defines the desired path decomposition.
\end{proof}

\clearpage
\begin{lem}
  \label{MinimalInducedConnectedSubgraph}\mbox{}\\*
  For every set $A$ of $\,k\geq2$ vertices in a connected graph $G$, every
  minimal induced connected subgraph $H$ of $\,G$ containing $A$ satisfies
  the following properties:
  \vspace*{-2.25ex}
  \begin{enumerate}[(1)]
    \setlength\itemsep{0ex}
  \item every (non-rooted) subtree of $\,H$ has at most $k$ leaves;
  \item $H$ has maximum degree at most $k$;
  \item $H$ has bandwidth (and hence pathwidth and treewidth) at most
    $k-1$;
  \item $H$ can be $2$-coloured with clustering $\bceil{\half k}$; and
  \item $H$ can be $2$-coloured with $\{\text{red},\text{blue}\}$ such that
    there are at most $k-2$ red vertices and the blue subgraph consists of
    at most $k-1$ pairwise disjoint paths.
  \end{enumerate}
\end{lem}

\preproof\begin{proof}
  Let $T$ be a spanning tree of $H$. By the minimality of $H$, every leaf
  of $T$ is in $A$. Thus~$T$ has at most $k$ leaves. Now let $S$ be any
  tree in~$H$. Extending~$S$ to a spanning tree of~$H$ cannot decrease the
  number of leaves, hence $S$ also has at most $k$ leaves.

  The closed neighbourhood of a vertex $v\in V(H)$ contains a tree with
  $\deg_H(v)$ leaves, proving $\deg_H(v)\leq k$.

  Let $T$ be a LexBFS spanning tree of $H$ rooted at a vertex $r$ in $A$.
  By the minimality of $H$, every leaf of $T$ is in $A$. Thus $T$ has at
  most $k-1$ leaves (the root does not count as a leaf). By
  \cref{LexBFSBandwidth}, $H$ has bandwidth, pathwidth and treewidth at
  most $k-1$.

  We now prove (4). We proceed by induction on $|V(H)|$. In the base case,
  $|V(H)|=|A|=k$ and the result is trivial. Now assume that $|V(H)|>k$.
  Thus $V(H-A)\ne\varnothing$, and by the minimality of $H$, every vertex
  in $H-A$ is a cut-vertex of $H$. Consider a leaf-block~$L$ of~$H$. Every
  vertex in $L$, except the one cut-vertex in $L$, is in $A$. There are at
  least two leaf-blocks. Thus $|V(L-v)|\leq\half k$ for some leaf block
  $L$, where $v$ is the one cut-vertex of $H$ in~$L$. Let $H'=H-V(L-v)$ and
  $A'=(A\setminus V(L))\cup\{v\}$. Then $H'$ is a minimal induced connected
  subgraph of $G$ containing $A'$, and $|A'|\leq k$. By induction, $H'$ has
  a $2$-colouring with clustering~$\bceil{\half k}$. Colour every vertex in
  $L\setminus\{v\}$ by the colour not assigned to $v$ in $H'$. Now $H$ is
  $2$-coloured with clustering $\bceil{\half k}$.

  It remains to prove (5). We proceed by induction on $k$. If $k=2$, then
  $H$ is a path between the two vertices in $A$. Colour every vertex in $H$
  blue, and we are done. So assume $k\geq 3$ and the result holds for
  $k-1$. Let $x$ be a vertex in $A$. By induction, every minimal induced
  connected subgraph $H'$ of $H$ containing $A\setminus\{x\}$ can be
  $2$-coloured with $\{\text{red},\text{blue}\}$ such that there are at
  most $k-3$ red vertices and the blue subgraph consists of at most $k-2$
  pairwise disjoint paths. If $x$ is in $H'$, then we are done. Otherwise,
  let $P$ be a shortest path between~$x$ and $H'$ in $H$. Say
  $P=x,\dots,u,v,w$, where $w$ is in $H'$. Then $v$ is the only vertex in
  $P-w$ adjacent to $H'$. Colour $v$ red, and colour $x,\dots,u$ blue. (It
  is possible that $x=v$, in which case $\{x,\dots,u\}=\varnothing$.) Then
  $\{x,\dots,u\}$ induces a path in $H$ that is not adjacent to~$H'$. By
  the minimality of $H$, we have $V(H)=V(H')\cup\{x,\dots,u,v\}$. Thus $H$
  is $2$-coloured with $\{\text{red},\text{blue}\}$ such that there are at
  most $k-2$ red vertices and the blue subgraph consists of at most $k-1$
  pairwise disjoint paths.
\end{proof}

We now prove the main result of this section. 

\begin{lem}\label{ConnectedSubgraph}\mbox{}\\*
  For every set $A$ of $\,k\geq2$ vertices in a connected graph $G$, there
  is an induced connected subgraph $H$ of $\,G$ containing $A$, such that:
  \vspace*{-2.25ex}
  \begin{enumerate}[(1)]
    \setlength\itemsep{0ex}
  \item $H$ has maximum degree at most $k$;
  \item $H$ has bandwidth (and hence pathwidth and treewidth) at most
    $k-1$;
  \item $H$ can be $2$-coloured with clustering $\bceil{\half k}$;
  \item $H$ can be $2$-coloured with $\{\text{red},\text{blue}\}$ such that
    there are at most $k-2$ red vertices and the blue subgraph consists of
    at most $k-1$ pairwise disjoint paths; and
  \item every vertex in $G$ has at most $2k-2$ neighbours in $V(H)$.
  \end{enumerate}
\end{lem}

\preproof\begin{proof}
  Let $T$ be a LexBFS spanning tree of $G$ rooted at some vertex $r\in A$.
  Let $S$ be the LexBFS subtree of $T$ consisting of all $ar$-paths in $T$,
  where $a\in A$. Every leaf of $S$ is in $A\setminus\{r\}$, implying that
  $S$ has at most $k-1$ leaves. By \cref{LexBFSDegree}, every vertex in $G$
  has at most $2k-2$ neighbours in $V(S)$. Let $H$ be a minimal induced
  connected subgraph of $G[V(S)]$ containing~$A$. The first four claims
  follow from \cref{MinimalInducedConnectedSubgraph}. Since
  $V(H)\subseteq V(S)$, \Cref{LexBFSDegree} means that every vertex in $G$
  has at most $2k-2$ neighbours in $V(H)$.
\end{proof}

\section{Decompositions of {\boldmath\texorpdfstring{$K_t$}{Kt}}-Minor-Free
  Graphs}
\label{ProofDecomposition}

\citet{HOQRS} introduced the following definition and proved the following
decomposition theorem. A connected partition $H_1,\dots,H_\ell$ has
\emph{width} $k$ if for each $i\in[\ell-1]$, each component of
$G-\bigl(V(H_1)\cup\dots\cup V(H_i)\bigr)$ is adjacent to at most $k$ of
the subgraphs $H_1,\dots,H_i$. Note that this implies that $H_{i+1}$ is
adjacent to at most $k$ of the subgraphs $H_1,\dots,H_i$ (since~$H_{i+1}$
is contained in some component of
$G-\bigl(V(H_1)\cup\dots\cup V(H_i)\bigr)$).

\begin{thm}[\citet{HOQRS}]\label{KtConnectedPartition}\mbox{}\\*
  Every $K_t$-minor-free graph $G$ has a connected partition
  $H_1,\dots,H_\ell$ with width $t-2$, such that each subgraph $H_i$ is
  induced by a BFS subtree of
  $\,G-\bigl(V(H_1)\cup\dots\cup V(H_{i-1})\bigr)$ with at most $t-3$
  leaves.
\end{thm}

The following similar decomposition theorem implies
\cref{BasicDecomposition}.

\clearpage
\begin{thm}\label{Decomposition}\mbox{}\\*
  For $t\geq4$, every $K_t$-minor-free graph $G$ has a connected partition
  $H_1,\dots,H_\ell$ with width $t-2$, such that for $i\in[\ell]$ the
  following holds.
  \vspace*{-2.25ex}
  \begin{enumerate}[(1)]
    \setlength\itemsep{0ex}
  \item The subgraph $H_i$ has the following properties:
    \vspace*{-1.25ex}
    \begin{enumerate}[(a)]
      \setlength\itemsep{0ex}
    \item $H_i$ has maximum degree at most $t-2$;
    \item $H_i$ has bandwidth, pathwidth and treewidth at most $t-3$;
    \item $H_i$ can be $2$-coloured with clustering $\bceil{\half(t-2)}$;
      and
    \item $H_i$ can be $2$-coloured with $\{\text{red},\text{blue}\}$ such
      that there are at most $t-4$ red vertices and the blue subgraph
      consists of at most $t-3$ pairwise disjoint paths.
    \end{enumerate}
    \vspace*{-1.25ex}
  \item Each component $C$ of $\,G-\bigl(V(H_1)\cup\dots\cup V(H_i)\bigr)$
    has the following properties.
    \vspace*{-1.25ex}
    \begin{enumerate}[(a)]
      \setlength\itemsep{0ex}
    \item At most $t-2$ subgraphs in $H_1,\dots,H_i$ are adjacent to $C$,
      and these subgraphs are pairwise adjacent. (This implies that at most
      $t-2$ subgraphs in $H_1,\dots,H_i$ are adjacent to $H_{i+1}$, and
      these subgraphs are pairwise adjacent.)
    \item Every vertex in $C$ is adjacent to at most $2t-6$ vertices in
      each of $\,H_1,\dots,H_i$. (This implies that every vertex in
      $H_{i+1}$ is adjacent to at most $2t-6$ vertices in each of
      $\,H_1,\dots,H_i$.)
    \end{enumerate}
  \end{enumerate}
\end{thm}

\preproof\begin{proof}
  We may assume that $G$ is connected. We construct $H_1,\dots,H_\ell$
  iteratively, maintaining properties (1) and (2). Let $H_1$ be the
  subgraph induced by a single vertex in $G$. Then (1) and~(2) hold for
  $i=1$.

  Assume that $H_1,\dots,H_i$ satisfy (1) and (2) for some $i\geq1$, but
  $V(H_1),\dots,V(H_i)$ do not partition $V(G)$. Let $C$ be a component of
  $G-\bigl(V(H_1)\cup\dots\cup V(H_i)\bigr)$. Let $Q_1,\dots,Q_k$ be the
  subgraphs in $H_1,\dots,H_i$ that are adjacent to $C$. By (2a),
  $Q_1,\dots,Q_k$ are pairwise adjacent and $k\leq t-2$. Since $G$ is
  connected, $k\geq1$.

  For $j\in[k]$, let $v_j$ be a vertex in $C$ adjacent to $Q_j$. If $k=1$,
  then let $H_{i+1}$ be the subgraph induced by $v_1$. It is easily seen
  that (1) is satisfied. For $k\geq 2$, by \cref{ConnectedSubgraph} with
  $k\leq t-2$, there is an induced connected subgraph $H_{i+1}$ of $C$
  containing $v_1,\dots,v_k$ that satisfies~(1).

  Consider a component $C'$ of
  $G-\bigl(V(H_1)\cup\dots\cup V(H_{i+1})\bigr)$. Either $C'$ is disjoint
  from $C$, or $C'$ is contained in $C$. If $C'$ is disjoint from $C$, then
  $C'$ is a component of $G-\bigl(V(H_1)\cup\dots\cup V(H_i)\bigr)$
  and~$C'$ is not adjacent to $H_{i+1}$, implying (2) is maintained for
  $C'$.

  Now assume $C'$ is contained in~$C$. Since every vertex in $C$ has at
  most $2t-6$ neighbours in each of $H_1,\dots,H_i$, every vertex in~$C'$
  has at most $2t-6$ neighbours in each of $H_1,\dots,H_i$. By
  \cref{ConnectedSubgraph}\,(5), every vertex in $C'$ also has at most
  $2t-6$ neighbours in $H_{i+1}$. Thus (2b) is maintained for~$C'$. The
  subgraphs in $H_1,\dots,H_{i+1}$ that are adjacent to $C'$ are a subset
  of $Q_1,\dots,Q_k,H_{i+1}$, which are pairwise adjacent. Suppose that
  $k=t-2$ and $C'$ is adjacent to all of $Q_1,\dots,Q_{t-2},H_{i+1}$. Then
  $C$ is adjacent to all of $Q_1,\dots,Q_{t-2}$. Contracting each of
  $Q_1,\dots,Q_{t-2},H_{i+1},C'$ into a single vertex gives $K_t$ as a
  minor of~$G$, a contradiction. Hence~$C'$ is adjacent to at most $t-2$ of
  $Q_1,\dots,Q_{t-2},H_{i+1}$, and property (2a) is maintained for~$C'$.
\end{proof}

Property (1d) in \cref{Decomposition}, along with a greedy
$(t-1)$-colouring of the subgraphs $H_1,\dots,H_\ell$, gives the following
results.

\begin{thm}\mbox{}\\*
  For $t\geq 4$, every $K_t$-minor-free graph has a $(2t-2)$-colouring such
  that for $t-1$ colours each monochromatic component has at most $t-4$
  vertices, and for the other $t-1$ colours each monochromatic component is
  a path.
\end{thm}

\begin{cor}\mbox{}\\*
  For $t\geq4$, every $K_t$-minor-free graph has a $(3t-3)$-colouring such
  that for $t-1$ colours, each monochromatic component has at most $t-4$
  vertices, and the other $2t-2$ colour classes are independent sets.
\end{cor}

The same greedy $(t-1)$-colouring of the subgraphs $H_1,\dots,H_\ell$,
together with \cref{Decomposition}\,(1b), gives the following result.

\begin{thm}\mbox{}\\*
  For $t\geq 4$, every $K_t$-minor-free graph has a $(t-1)$-colouring such
  that each monochromatic component has treewidth at most $t-3$.
\end{thm}

Note that \citet{Detal04} proved that for every proper minor-closed class
of graphs, every graph in that class has a \mbox{$2$-colouring} such that
each monochromatic component has bounded treewidth. Their proof again uses
the Graph Minor Structure Theorem, leading to a very large bound on the
treewidth.

Property (2a) in \cref{Decomposition} means that if $Q$ is the graph
obtained $G$ by contracting each subgraph $H_i$ to a single vertex, then
$Q$ is chordal with no $K_t$-subgraph, and thus with treewidth at most
$t-2$. Indeed, $H_1,\dots,H_\ell$ defines an elimination ordering of $Q$.
In the language of \citet{ReedSeymour-JCTB98}, $H_1,\dots,H_\ell$ is a
\emph{chordal decomposition} with \emph{touching pattern} $Q$. We only need
that $Q$ is $(t-2)$-degenerate for \cref{MonoComponents,MaxDegree}, but it
is interesting that, in fact,~$Q$ has treewidth at most $t-2$.

Even though we do not use it explicitly in this paper, it is an interesting
aspect of our decomposition that the superstructure (that is, $Q$) has
bounded treewidth, as does each piece of the decomposition. There are
several other properties in \Cref{Decomposition} we do not use, but we
mention them since they might be useful for other applications.

\section{Excluding a Complete Bipartite Minor}
\label{K2tK3t}

This section presents decomposition results analogous to
\cref{Decomposition} for $K_{s,t}$-minor-free graphs, leading to bounds on
the defective and clustered chromatic number. Those decomposition results
and the more technical proofs can be found towards the end of the section.

In fact, for most of this section we will consider the larger classes of
$K^*_{s,t}$-minor-free graphs, where $K^*_{s,t}$ is the complete join of
$K_s$ and $\overline{K_t}$. We start with $s\in\{1,2,3\}$ before
considering the general case. Graphs with no $K_{1,t}$ minor (note that
$K_{1,t}=K^*_{1,t}$) are easily coloured. Every such graph has maximum
degree at most $t-1$, and is therefore $1$-colourable with defect $t-1$.
Moreover, every BFS layer has at most $t-1$ vertices, so alternately
colouring the BFS layers gives a $2$-colouring with clustering $t-1$.

Next consider the $s=2$ case. \citet{OOW16} proved that every
$K^*_{2,t}$-minor-free graph is $2$-colourable with defect $\Oh{t^3}$. Our
decomposition results imply the following improvement.

\begin{thm}\label{K2tColouring}\mbox{}\\*
  Every $K^*_{2,t}$-minor-free graph is $2$-colourable with defect $2t-2$.
\end{thm}

The decomposition results for $K^*_{2,t}$-minor-free graphs also imply that
every is $4$-colourable with clustering $t-1$. This result can be improved
as follows. The proof is inspired by a method of \citet{Gon11}.

\begin{thm}\label{thmGon}\mbox{}\\*
  Every $K^*_{2,t}$-minor-free graph $G$ is $3$-colourable with clustering
  $t-1$. Moreover, for each edge $vw$ of $\,G$, there is such a
  $3$-colouring in which $v$ and $w$ are both isolated in their respective
  monochromatic subgraphs.
\end{thm}

Now consider $K^*_{3,t}$-minor-free graphs. \citet{OOW16} proved that the
defective chromatic number of $K^*_{3,t}$-minor-free graphs equals 3. In
particular, every $K^*_{3,t}$-minor-free graph is $3$-colourable with
defect $\Oh{t^4}$. Our decomposition results again imply an improvement.

\begin{thm}\label{K3tColouring}\mbox{}\\*
  Every $K^*_{3,t}$-minor-free graph is $3$-colourable with defect $4t$,
  and is $6$-colourable with clustering~$2t$.
\end{thm}

It follows from Euler's Formula that graphs with Euler genus $g$ exclude
$K_{3,2g+3}$ as a minor. Thus the second part of \cref{K3tColouring} is
related to the results of \citet{EO16} and \citet{KT12} that every graph of
Euler genus~$g$ can be \mbox{$5$-coloured} with clustering $\Oh{g}$.
\citet{KMRV97} constructed planar graphs that cannot be $3$-coloured with
bounded clustering. We conjecture that every $K_{3,t}$-minor-free graph is
$4$-colourable with clustering $f(t)$, for some function~$f$.

It is possible to improve the bound on the cluster size for the
$6$-colouring result in \cref{K3tColouring}. In a
\mbox{$K^*_{3,t}$-minor}-free graph, every BFS layer induces a
$K^*_{2,t}$-minor-free graph, which is $3$-colourable with clustering $t-1$
by \cref{thmGon}. Using disjoint sets of three colours for alternate BFS
layers gives a $6$-colouring with clustering $t-1$.

Finally, in this section we consider general $K_{s,t}$-minor-free graphs.
\citet{OOW16} proved that the defective chromatic number of
$K_{s,t}$-minor-free graphs equals $s$. We show that the clustered
chromatic number of $K_{s,t}$-minor-free graphs is at least $s+1$, thus
generalising the above-mentioned lower bound of \citet{KMRV97}.

\begin{prop}\label{KstLowerBound}\mbox{}\\*
  For every $s\geq1$, $t\geq\max\{s,3\}$ and $c\geq1$, there is a
  $K_{s,t}$-minor-free graph $G_s$ such that every $s$-colouring of $\,G_s$
  has a monochromatic component of order greater than $c$.
\end{prop}

\preproof\begin{proof}
  Define $G_s$ recursively as follows. Let $G_1$ be the path on $c+1$
  vertices. For $s\geq2$, let~$G_s$ be the graph obtained from $c$ disjoint
  copies of $G_{s-1}$ by adding one dominant vertex.

  We claim that $G_s$ is not $s$-colourable with clustering $c$. We prove
  this claim by induction on $s\geq1$. Obviously, $G_1$ is not
  $1$-colourable with clustering $c$. Now assume that $s\geq2$
  and~$G_{s-1}$ is not $(s-1)$-colourable with clustering~$c$. Suppose that
  $G_{s}$ has an $s$-colouring with clustering~$c$. Say the dominant vertex
  in $G_s$ is coloured black. At most $c-1$ copies of~$G_{s-1}$ contain a
  black vertex, which implies that at least one copy has no black vertex.
  Thus~$G_{s-1}$ has an $(s-1)$-colouring with clustering $c$, which is a
  contradiction. Hence $G_{s}$ is not $s$-colourable with clustering $c$,
  as claimed.

  It remains to show that $G_{s}$ is $K_{s,t}$-minor-free with
  $t\geq\max\{s,3\}$. We do so by induction on $s\geq1$. $G_1$ is a path,
  and therefore contains no $K_{1,3}$ minor. $G_2$ is outerplanar, and
  therefore contains no $K_{2,3}$ minor. $G_3$ is planar, and therefore
  contains no $K_{3,3}$ minor.

  Now assume that $s\geq4$ and $G_{s-1}$ contains no $K_{s-1,s-1}$ minor,
  but $G_{s}$ contains a $K_{s,s}$ minor. Let $v$ be the dominant vertex in
  $G_s$. We may assume that $v$ is the entire image of one vertex in
  the~$K_{s,s}$ minor in~$G$. Since~$K_{s,s}$ is $2$-connected, the
  $K_{s,s}$ minor is contained in one copy of~$G_{s-1}$ plus~$v$. Deleting
  any one vertex from $K_{s,s}$ gives a subgraph that contains a
  $K_{s-1,s-1}$ subgraph. Thus~$G_{s-1}$ contains a $K_{s-1,s-1}$ minor,
  which is a contradiction. We conclude that for $s\geq4$, $G_s$ has
  no~$K_{s,s}$ minor, so certainly no $K_{s,t}$ minor with $t\geq s$
  ($=\max\{s,3\}$).
\end{proof}

Determining the clustered chromatic number of $K_{s,t}$-minor-free graphs
is an open problem. \cref{KstLowerBound} provides a lower bound of $s+1$.
Since $K_{s,t}$-minor-free graphs are defectively $s$-colourable
\citep{OOW16}, \cref{Defective2Fragmented} implies an upper bound of $3s$.
In general, for every graph $H$, it is possible that the clustered
chromatic number of $H$-minor-free graphs is at most one more than the
defective chromatic number of $H$-minor-free graphs.

We now give the structural results and proofs of the above statements in
this section. All the results in this section are based on LexBFS, so we
first present the following general lemma. Recall the definition of the
width of a connected partition from the beginning of
\Cref{ProofDecomposition}.

\clearpage
\begin{lem}\label{PartitionColouring}\mbox{}\\*
  Suppose that a graph $G$ has a connected partition $H_1,\dots,H_\ell$
  with width $k$. If each subgraph~$H_i$ is induced by a BFS subtree of
  $\,G-V\bigl(V(H_1)\cup\dots\cup V(H_{i-1})\bigr)$ with at most $p$
  leaves, then $G$ is $(k+1)$-colourable with defect $3p-1$, and $G$ is
  $(2k+2)$-colourable with clustering $p$.

  If, in addition, each subgraph $H_i$ is induced by a LexBFS subtree of
  $\,G-\bigl(V(H_1)\cup\dots\cup V(H_{i-1})\bigr)$ with at most $p$ leaves,
  then $G$ is $(k+1)$-colourable with defect $2p$.
\end{lem}

\preproof\begin{proof}
  Colour the subgraphs $H_1,\dots,H_\ell$ greedily in this order, such that
  adjacent subgraphs receive distinct colours. Since the partition has
  width~$k$, $k+1$ colours suffice. Colour each vertex in $H_i$ by the
  colour assigned to $H_i$. In each subgraph $H_i$ each BFS layer has at
  most~$p$ vertices. Since a vertex in a BFS subtree has neighbours in its
  own layer and in the two layers below and above its own layer only, $H_i$
  has maximum degree at most $3p-1$. Hence $G$ is $(k+1)$-colourable with
  defect $3p-1$. Moreover, if each subgraph~$H_i$ is induced by a LexBFS
  subtree of $G-\bigl(V(H_1)\cup\dots\cup V(H_{i-1})\bigr)$ with at most
  $p$ leaves, then by \cref{LexBFSDegree}, $H_i$ has maximum degree $2p$.

  For the clustering claim, alternately $2$-colour the BFS layers in each
  $H_i$, and take the product with the $(k+1)$-colouring of
  $H_1,\dots,H_\ell$ to produce a $(2k+2)$-colouring of $G$ with
  clustering~$p$.
\end{proof}

As an aside, note that \citet{HOQRS} proved that every planar graph has a
connected partition $H_1,\dots,H_n$ with width 2, such that each subgraph
$H_i$ is a shortest path in $G-\bigl(V(H_1)\cup\dots\cup V(H_{i-1})\bigr)$.
Thus, \cref{PartitionColouring} with $k=2$ and $p=1$ implies that planar
graphs are $3$-colourable with defect 2, which is the best possible result
for defective $3$-colouring of planar graphs, first proved by
\citet{CCW86}. In fact, each monochromatic component is a path, which was
previously proved by \citet{Goddard91} and \citet{Poh90}.

For $K^*_{2,t}$-minor-free graphs we have the following.

\begin{lem}\label{K2tPartition}\mbox{}\\*
  Every $K^*_{2,t}$-minor-free graph $G$ has a connected partition
  $H_1,\dots,H_\ell$ with width $1$, such that each subgraph $H_i$ is
  induced by a LexBFS subtree of $\,G-\bigl(V(H_1)\cup\dots\cup
  V(H_{i-1})\bigr)$ with at most $t-1$ leaves.
\end{lem}

\preproof\begin{proof}
  We may assume that $G$ is connected. We construct $H_1,\dots,H_\ell$
  iteratively. Let $H_1$ be the subgraph induced by a single vertex in $G$.

  Assume that $H_1,\dots,H_i$ are defined for some $i\geq1$, and $C$ is a
  component of $G-\bigl(V(H_1)\cup\dots\cup V(H_i)\bigr)$ adjacent to one
  of $H_1,\dots,H_i$. (Since $G$ is connected, $C$ is adjacent to at least
  one of those subgraphs.) So $C$ is adjacent to $H_a$, for some $a\in[i]$,
  and to no other subgraph in $H_1,\dots,H_i$; let~$A$ be the set of
  vertices in $C$ adjacent to $H_a$, and let $r$ be a vertex in~$A$.
  Let~$S$ be a LexBFS subtree of $C$ rooted at $r$, such that every vertex
  in $A$ is in $S$, and subject to this property, $|V(S)|$ is minimal. Thus
  every leaf of $S$ is in $A$. Let $S_0$ be the subtree of~$S$ obtained by
  deleting the leaves. If $S$ has at least $t$ leaves, then a $K^*_{2,t}$
  minor is obtained by contracting~$H_a$ to a vertex and contracting $S_0$
  to a vertex. Thus $S$ has at most $t-1$ leaves. Let~$H_{i+1}$ be the
  subgraph of $C$ induced by $V(S)$. Since every vertex in $A$ is in~$S$,
  every component of $G-\bigl(V(H_1)\cup\dots\cup V(H_{i+1})\bigr)$ is
  adjacent to at most one of $H_1,\dots,H_{i+1}$. Iterating this process
  gives the desired partition.
\end{proof}

\cref{PartitionColouring,K2tPartition} immediately imply
\Cref{K2tColouring} and show that $K^*_{2,t}$-minor-free graphs are
$4$-colourable with clustering $t-1$. As expressed in \Cref{thmGon}, this
can be improved to a $3$-colouring with the same clustering bounds, as we
now prove.

\begin{proof}[Proof of \Cref{thmGon}.]
  We proceed by induction on $|V(G)|$. The claim is trivial if $|V(G)|\leq
  t+1$. Now assume that $vw$ is an edge in a $K^*_{2,t}$-minor-free graph
  $G$, and the result holds for \mbox{$K^*_{2,t}$-minor}-free graphs with
  fewer vertices than $G$. If $\deg_G(v)=1$, then by induction $G-v$ has a
  $3$-colouring in which $w$ is isolated in its monochromatic subgraph.
  Assign $v$ a colour not assigned to $w$. We obtain the desired colouring
  of $G$.

  Now assume that $\deg_G(v)\geq2$ and, similarly, $\deg_G(w)\geq2$. Let
  $A$ and $B$ be disjoint sets of vertices in $G$ such that $v\in A$ and
  $w\in B$, $G[A]$ and $G[B]$ are connected, and $vw$ is the only edge
  between $A$ and $B$, and subject to these properties, $|A\cup B|$ is
  maximum. The sets $A$ and~$B$ are well-defined, since $A=\{v\}$ and
  $B=\{w\}$ satisfy the conditions. Let $Z$ be the set of vertices in
  $V(G)\setminus(A\cup B)$ adjacent to both $A$ and $B$, and let
  $Y=V(G)\setminus(A\cup B\cup Z)$.

  If $|Z|\geq t$, then contracting $A$ and $B$ into single vertices gives a
  $K^*_{2,t}$ minor. Thus $|Z|\leq t-1$. Since $G[A]$ is connected and
  every vertex in $Z$ is adjacent to $A$, $G[A\cup Z]$ is connected.
  Similarly, $G[B\cup Z]$ is connected.

  Let $G_1$ be obtained from $G$ by contracting $G[B\cup Z]$ into a single
  vertex $x$. Note that $vx$ is an edge of $G_1$. Let $G_2$ be obtained
  from $G$ by contracting $G[A\cup Z]$ into a single vertex $y$. Note that
  $wy$ is an edge of $G_2$. Since $G_1$ and $G_2$ are minors of $G$, they
  both contain no $K^*_{2,t}$ minor. Since $\deg_G(v)\geq2$ and
  $\deg_G(w)\geq2$, both $G_1$ and $G_2$ have fewer vertices than $G$.

  By induction, $G_1$ is $3$-colourable with clustering $t-1$ such that $v$
  and $x$ are both isolated in their respective monochromatic subgraphs,
  and $G_2$ is $3$-colourable with clustering $t-1$ such that $w$ and $y$
  are both isolated in their respective monochromatic subgraphs. Permute
  the colours in $G_2$ so that $x\in V(G_1)$ and $y\in V(G_2)$ receive the
  same colour, and $v\in V(G_1)$ and $w\in V(G_2)$ receive distinct
  colours.

  Let $G_3$ be obtained from $G$ by contracting $G[A\cup B\cup Z]$ into a
  single vertex $z$. Note that $V(G_3)=Y\cup\{z\}$. By induction, $G_3$ is
  $3$-colourable with clustering $t-1$ such that $z$ is isolated in its
  colour class. Permute the colours in $G_3$ so that $z$ receives the same
  colour as $x\in V(G_1)$, which is the same colour assigned to
  $y\in V(G_2)$.

  Colour each vertex in $Z$ by the colour assigned to $x$ and $y$. Colour
  each vertex in $A$ by its colour in $G_1$. Colour each vertex in $B$ by
  its colour in $G_2$. Finally, colour each vertex in $Y$ by its colour in
  $G_3$.
  
  Since $x$ is isolated in its monochromatic subgraph in $G_1$, $y$ is
  isolated in its monochromatic subgraph in $G_2$, and $z$ is isolated in
  its monochromatic subgraph in $G_3$, every monochromatic component
  intersecting $Z$ is contained in $Z$, and thus has at most $t-1$
  vertices. Since $vw$ is the only edge between $A$ and $B$, and $v$ and
  $w$ are assigned distinct colours, every monochromatic component that
  intersects $A$ is contained in $A$, and therefore by induction has at
  most $t-1$ vertices. Similarly, every monochromatic component that
  intersects $B$ is contained in $B$, and therefore by induction has at
  most $t-1$ vertices.
\end{proof}

The following lemma is used in our decomposition result for
$K^*_{3,t}$-minor-free graphs.

\begin{lem}\label{AB}\mbox{}\\*
  For every connected graph $G$, non-empty sets $A,B\subseteq V(G)$, and
  integer $t\geq1$, 
  \vspace*{-2.25ex}
  \begin{enumerate}[(1)]
    \setlength\itemsep{0ex}
  \item $G$ has a LexBFS subtree $T$ with at most $2t+1$ leaves, such that
    $T$ intersects both $A$ and~$B$, and $V(T)$ separates $A$ and $B$; or
  \item $G$ has a $K_{1,t}$ minor with every branch set intersecting both
    $A$ and $B$.
  \end{enumerate}
\end{lem}

\preproof\begin{proof}
  Let $r$ be a vertex in $A$. Let $X$ be a LexBFS spanning tree of $G$
  rooted at $r$. For a set $L\subseteq V(G)$, let $T_L$ be the subtree of
  $X$ consisting of the union of all paths in $X$ between $L$ and $r$.
  Choose $L\subseteq V(G)$ so that $V(T_L)$ is an $AB$-separator, and
  subject to this property, $|L|$ is minimum. This is well-defined since if
  $L=V(G)$, then $V(T_L)=V(G)$. By the minimality of~$|L|$, every vertex in
  $L$ is a leaf of $T_L$. And by the definition of $T_L$, every leaf of
  $T_L$ is in $L$.

  For each $x\in L$, let $p_x$ be the vertex closest to $x$ in $T_L$, such
  that $\deg_T(p_x)\geq3$ or $p_x=r$. Let $Q_x$ be the path in $T_L$
  between $x$ and $p_x$ not including $p_x$. We call $Q_x$ the \emph{leaf
    path at $x$}. Let $T_0=T_L-\bigcup_{x\in L}V(Q_x)$. Let $H$ be the
  graph obtained from $G-V(T_0)$ by contracting the leaf path $Q_x$
  corresponding to each $x\in L$ into a single vertex $y_x$. 
  We consider $A$ to also be a set of vertices in $H$, where a vertex $y_x$ is in $A$ if any vertex of $Q_x$ is in $A$, and similarly for $B$.  Let~$S$ be a  minimum $AB$-separator in $H$.

  First suppose that $|S|\geq t+1$. By Menger's Theorem, there are $t+1$
  pairwise disjoint \mbox{$AB$-paths} $Z_1,\dots,Z_{t+1}$ in $H$. Since
  $V(T_L)$ is an $AB$-separator in $G$, each $Z_i$ contains $y_x$ for some
  $x\in L$, and each vertex $y_x$ is on at most one path $Z_i$. For
  $i\in[t]$, if $Z_i$ contains $y_x$, then contract $Z_i\cup Q_x$ into a
  single vertex. If $Z_{t+1}$ contains $y_x$, then contract
  $Z_{t+1}\cup Q_x\cup T_0$ into a single vertex. We obtain a $K_{1,t}$
  minor with every branch set intersecting both $A$ and $B$ (since
  each~$Q_x$ is adjacent to $T_0$), and (2) holds. 

%
%
%

Now assume that $|S|\leq t$. 
 Let~$S_1$ be the set of vertices $x\in L$  such that $y_x$ is in $S$. 
  Let~$S_2$ be the set of vertices  in $G-S_1$ that correspond to vertices in $S$. Thus $|S|=|S_1|+|S_2|$. 
 Let $Z$ be the set of vertices $z\in V(T_0)$ such that $z=p_x$ for some $x\in L\setminus S_1$, and $z\neq p_x$ for all $x\in S_1$. 
 Let~$T'=T_{L'}$,  where  $L'=S_1\cup S_2\cup Z$. 
  Since~$S$ separates $A$ and $B$ in
  $H$, and $T'$ contains $T_0\cup S_1\cup S_2$ along with $Q_x$ for each $x\in S_1$,
  it follows that~$V(T')$ separates $A$ and $B$ in $G$. 

By the definition of  $p_x$, for each vertex $z\in Z$ there are at least two vertices $x$ and $x'$ in $L\setminus S_1$ for which $z=p_x=p_{x'}$. Thus $|L\setminus S_1| \geq 2|Z|$. By the choice of $L$, we can argue
  \[|L|\leq |L'| = |S_1|+|S_2|+|Z| \leq   |S|+ \half|L\setminus S_1| \leq t+\half |L|.\]
  Hence $|L|\leq2t$. Thus $T_L$ is a LexBFS subtree with at most $2t$
  leaves, such that $V(T_L)$ separates~$A$ and $B$. Let $T$ be obtained
  from $T_L$ by adding a shortest path in $X$ from $T_L$ to $B$. Then $T$
  is a LexBFS subtree with at most $2t+1$ leaves, such that $T$ intersects
  both $A$ and $B$, and $V(T)$ separates $A$ and $B$.
\end{proof}

We are now ready to prove the following structural lemma.

\begin{lem}\label{K3tPartition}\mbox{}\\*
  Every $K^*_{3,t}$-minor-free graph $G$ has a connected partition
  $H_1,\dots,H_\ell$ with width $2$, such that each subgraph $H_i$ is
  induced by a LexBFS subtree of
  $\,G-\bigl(V(H_1)\cup\dots\cup V(H_{i-1})\bigr)$ with at most $2t+1$
  leaves.
\end{lem}

\preproof\begin{proof}
  We again may assume that $G$ is connected. We construct
  $H_1,\dots,H_\ell$ iteratively, maintaining the property that for each
  $i\in[\ell-1]$, each component $C$ of
  $G-\bigl(V(H_1)\cup\dots\cup V(H_i)\bigr)$ is adjacent to at most two of
  $H_1,\dots,H_i$, and if $C$ is adjacent to $H_a$ and $H_b$, for some
  distinct $a,b\in[i]$, then $H_a$ and~$H_b$ are adjacent. Call this
  property $(\star)$.

  Assume that $H_1,\dots,H_i$ is defined, and $C$ is a component of
  $G-\bigl(V(H_1)\cup\dots\cup V(H_i)\bigr)$. (Hence~$C$ satisfies property
  $(\star)$.)

  Suppose $C$ is adjacent to $H_a$, for some $a\in[i]$, and to no other
  subgraph in $H_1,\dots,H_i$. Let~$H_{i+1}$ be a subgraph of $C$ induced
  by one vertex adjacent to $H_a$. Let $C'$ be a component of
  $G-\bigl(V(H_1)\cup\dots\cup V(H_{i+1})\bigr)$. If $C'$ is a component of
  $G-\bigl(V(H_1)\cup\dots\cup V(H_i)\bigr)$, then $C'$ is not adjacent to
  $H_i$, and $(\star)$ is maintained for $C'$. Otherwise $C'$ is a
  component of $C-V(H_{i+1})$, and $C'$ is adjacent to $H_{i+1}$ and
  possibly $H_a$. Since $H_{i+1}$ and $H_a$ are adjacent, $(\star)$ holds
  for~$C'$.

  Now assume that $C$ is adjacent to $H_a$ and $H_b$, for some distinct
  $a,b\in[i]$, and to no other subgraph in $H_1,\dots,H_i$. Let $A$ be the
  set of vertices in $C$ adjacent to $H_a$, and let $B$ be the set of
  vertices in $C$ adjacent to $H_b$. By \cref{AB} above we have: (1)~$C$
  has a LexBFS subtree~$T$ separating $A$ and $B$, such that $T$ intersects
  both $A$ and $B$, and $T$ has at most $2t+1$ leaves, or (2)~$C$ has a
  $K_{1,t}$ minor with every branch set intersecting both $A$ and $B$. In
  case (1), let $H_{i+1}$ be the subgraph of $C$ induced by $V(T)$. Since
  $T$ intersects both $A$ and $B$, the subgraph~$H_{i+1}$ is adjacent to
  both $H_a$ and $H_b$. Let $C'$ be a component of
  $G-\bigl(V(H_1)\cup\dots\cup V(H_{i+1})\bigr)$. If~$C'$ is a component of
  $G-\bigl(V(H_1)\cup\dots\cup V(H_i)\bigr)$, then $C'$ is not adjacent to
  $H_{i+1}$, and $(\star)$ is maintained for $C'$. Otherwise, $C'$ is a
  component of $C-V(H_{i+1})$. Then $C'$ is adjacent to~$H_{i+1}$ and at
  most one of $H_a$ and $H_b$ (since $V(T)$ separates $A$ and $B$). Thus
  property $(\star)$ holds for~$C'$ (since $H_{i+1}$ is adjacent to both
  $H_a$ and $H_b$).

  In Case (2), with $H_a$ and $H_b$ we obtain a $K^*_{3,t}$ minor in $G$,
  which is a contradiction.
\end{proof}

\section{Generalised Colouring Numbers}
\label{GeneralisedColouring}

This section presents bounds on generalised colouring numbers, first
introduced by \citet{KY03}. Generalised colouring numbers are important
because they characterise bounded expansion classes \citep{Zhu09}, they
characterise nowhere dense classes \citep{GKRSS15}, and have several
algorithmic applications such as the constant-factor approximation
algorithm for domination number by \citet{Dvorak13}, and the almost
linear-time model-checking algorithm of Grohe et al.~\citep{GKS14}. They
also interpolate between degeneracy and treewidth (strong colouring
numbers) and between degeneracy and treedepth (weak colouring numbers). See
\citep{HOQRS,KPRS16,Sparsity} for more details.

For a graph $G$, linear ordering $\preceq$ of $V(G)$, vertex $v\in V(G)$,
and integer $r\geq1$, let $\sreach_r(G,\preceq,v)$ be the set of vertices
$x\in V(G)$ for which there is a path $v=w_0,w_1,\dots,w_{r'}=x$ of length
$r'\in[0,r]$ such that $x\preceq v$ and $v\prec w_i$ for all $i\in[r-1]$.
Similarly, let $\wreach_r(G,\preceq,v)$ be the set of vertices $x\in V(G)$
for which there is a path $v=w_0,w_1,\dots,w_{r'}=x$ of length $r'\in[0,r]$
such that $x\preceq v$ and $x\prec w_i$ for all $i\in[r'-1]$. For a graph
$G$ and integer $r\geq1$, the \emph{$r$-strong colouring number}
$\scol_r(G)$ of~$G$ is the minimum integer $k$ such that there is a linear
ordering~$\preceq$ of $V(G)$ with $|\sreach_r(G,\preceq,v)|\leq k$ for each
vertex $v$ of $G$. Similarly, the \emph{$r$-weak colouring number}
$\wcol_r(G)$ is the minimum integer $k$ such that there is a linear
ordering $\preceq$ of~$V(G)$ with $|\wreach_r(G,\preceq,v)|\leq k$ for each
vertex $v$ of $G$.

The following lemma is implicitly proved by \citet{HOQRS}.

\begin{lem}[\citet{HOQRS}]\label{BFScolouring}\mbox{}\\*
  Let $H_1,\dots,H_\ell$ be a connected partition of a graph $G$ with width
  $k$, such that there exists $p$ such that for $i\in[\ell]$,
  $V(H_i)=V(P_{i,1})\cup\dots\cup V(P_{i,p_i})$, where $p_i\leq p$ and each
  $P_{i,j}$ is a shortest path in
  $G-\Bigl(\bigl(V(H_1)\cup\dots\cup V(H_{i-1})\bigr)\cup
  \bigl(V(P_{i,1})\cup\dots\cup V(P_{i,j-1})\bigr)\Bigr)$. Then the
  generalised colouring numbers of~$\,G$ satisfy for every $r\geq1$:
  \[\scol_r(G)\leq p(k + 1)(2r+1)\qquad \text{and}\qquad
  \wcol_r(G)\leq p\binom{r+k}{k}(2r+1).\]
\end{lem}

Note that the conditions on the paths $P_{i,j}$ in the lemma are implied if
$H_i$ is induced by a BFS subtree with at most $p$ leaves in
$G-\bigl(V(H_1)\cup\dots\cup V(H_{i-1})\bigr)$.

For example, combining \Cref{BFScolouring} with a variant of
\cref{KtConnectedPartition}, \citet{HOQRS} proved that every
$K_t$-minor-free graph $G$ satisfies:
\[\scol_r(G)\leq \binom{t-1}{2}(2r+1)\qquad \text{and}\qquad
\wcol_r(G)\leq (t-3)\binom{r+t-2}{t-2}(2r+1).\]

\cref{BFScolouring,K2tPartition} imply:

\begin{thm}\label{K2t}\mbox{}\\*
  For every $K^*_{2,t}$-minor-free graph $G$ and every $r\geq1$,
  \[\scol_r(G)\leq 2(t-1)(2r+1)\qquad \text{and}\qquad
  \wcol_r(G)\leq (t-1)(r+1)(2r+1).\]
\end{thm}

And \cref{BFScolouring,K3tPartition} imply:

\begin{thm}\label{K3t}\mbox{}\\*
  For every $K^*_{3,t}$-minor-free graph $G$ and every $r\geq1$,
  \[\scol_r(G)\leq 3(2t+1)(2r+1)\qquad \text{and}\qquad
  \wcol_r(G)\leq (2t+1)\binom{r+2}{2}(2r+1).\]
\end{thm}

Since graphs with Euler genus $g$ exclude $K_{3,2g+3}$ as a minor,
\cref{K3t} implies that for every graph $G$ with Euler genus~$g$,
\[\scol_r(G)\leq 3(4g+7)(2r+1)\qquad \text{and}\qquad
\wcol_r(G)\leq (4g+7)\binom{r+2}{2}(2r+1).\]
These result are within a constant factor of the best known bounds for
graphs of Euler genus~$g$, due to \citet{HOQRS}. Note that \cref{K3t}
applies to a broader class of graphs than those with bounded Euler genus.
For example, the disjoint union of $g+1$ copies of $K_5$ has Euler genus
$g+1$, but contains no $K_{3,3}$ minor. It is easy to construct
$3$-connected examples as well.

We conjecture that \cref{K2t,K3t} can be generalised as follows:

\begin{conj}\label{KstConjectureColouring}\mbox{}\\*
  There exists a function $f$ such that for every $K^*_{s,t}$-minor-free
  graph $G$ and every $r\geq1$,
  \[\wcol_r(G)\leq f(s,t)\,r^s.\]
\end{conj}

\cref{KstConjectureColouring} would be implied by \cref{BFScolouring} and
the following conjecture.

\begin{conj}\label{KstConjectureDecomposition}\mbox{}\\*
  For all $t\geq s\geq1$, there exists an integer $p$, such that every
  $K^*_{s,t}$-minor-free graph $G$ has a connected partition
  $H_1,\dots,H_\ell$ with width $s-1$, such that for $i\in[\ell]$,
  $V(H_i)=V(P_{i,1})\cup\dots\cup V(P_{i,p_i})$, where $p_i\leq p$ and each
  $P_{i,j}$ is a shortest path in
  $G-\Bigl(\bigl(V(H_1)\cup\dots\cup V(H_{i-1})\bigr)\cup
  \bigl(V(P_{i,1})\cup\dots\cup V(P_{i,j-1})\bigr)\Bigr)$.
\end{conj}

We now show that \cref{KstConjectureColouring} is true with $r^s$ replaced
by $r^{s+1}$.

\begin{prop}\label{KstColouring}\mbox{}\\*
  For every $K^*_{s,t}$-minor-free graph $G$ and every $r\geq1$, we have
  \[\scol_r(G)\leq s(s+1)(t-1)(2r+1)\quad \text{and}\quad
  \wcol_r(G)\leq s(t-1)\binom{r+s}{s}(2r+1).\]
\end{prop}

\cref{KstColouring} follows from \cref{BFScolouring} and the next lemma.

\begin{lem}\label{LemKst}\mbox{}\\*
  Every $K^*_{s,t}$-minor-free graph has a connected partition
  $H_1,\dots,H_\ell$ with width $s$, such that for $i\in[\ell]$,
  $V(H_i)=V(P_{i,1})\cup\dots\cup V(P_{i,p_i})$, where $p_i\leq s(t-1)$ and
  each $P_{i,j}$ is a shortest path in
  $G-\Bigl(\bigl(V(H_1)\cup\dots\cup V(H_{i-1})\bigr)\cup
  \bigl(V(P_{i,1})\cup\dots\cup V(P_{i,j-1})\bigr)\Bigr)$.
\end{lem}

\preproof\begin{proof}
  Once more we may assume that $G$ is connected. We construct
  $H_1,\dots,H_\ell$, maintaining the property that for each
  $i\in[\ell-1]$, each component $C$ of
  $G-\bigl(V(H_1)\cup\dots\cup V(H_i)\bigr)$ is adjacent to at most $s$
  subgraphs in $H_1,\dots,H_i$, and that the subgraphs~$C$ is adjacent to
  are also pairwise adjacent. Call this property $(\star)$.

  Assume that $H_1,\dots,H_i$ is defined, and $C$ is a component of
  $G-\bigl(V(H_1)\cup\dots\cup V(H_i)\bigr)$. (Hence~$C$ satisfies property
  $(\star)$.) Let $Q_1,\dots,Q_k$ be the subgraphs in $H_1,\dots,H_i$ that
  are adjacent to~$C$. Thus $Q_1,\dots,Q_k$ are pairwise adjacent and
  $k\leq s$.

  Since $G$ is connected, $k\geq1$. For $j\in[k]$, let $A_j$ be the set of
  vertices in $C$ adjacent to~$Q_j$. Each~$A_j$ is non-empty. Let
  $\{F_1,\dots,F_m\}$ be a maximal set of pairwise disjoint connected
  subgraphs constructed as follows. The subgraph~$F_1$ is induced by a
  minimal BFS subtree $S_1$ in~$C$ rooted at some vertex $v\in V(C)$ and
  with~$S_1$ intersecting all of $A_1,\dots,A_k$. For $j\geq1$, $F_{j+1}$
  is induced by a minimal BFS subtree $S_{j+1}$ in
  $C-\bigl(V(F_1)\cup\dots\cup V(F_j)\bigr)$ rooted at some vertex~$v$ that
  is adjacent to $F_1\cup\dots\cup F_j$, and with $S_{j+1}$ intersecting
  all of $A_1,\dots,A_k$. By minimality, each $S_j$ has at most $k\leq s$
  leaves. Thus each~$S_j$ is the union of at most $s$ shortest paths in
  $C-\bigl(V(F_1)\cup\dots\cup V(F_{j-1})\bigr)$.

  Suppose that $k\leq s-1$. Let $H_{i+1}=F_1$. Then $H_{i+1}$ satisfies the
  claim. Consider a component~$C'$ of
  $G-\bigl(V(H_1)\cup\dots\cup V(H_{i+1})\bigr)$. If $C'$ is disjoint from
  $C$, then $C'$ is a component of
  $G-\bigl(V(H_1)\cup\dots\cup V(H_i)\bigr)$ and $C'$ is not adjacent to
  $H_{i+1}$. Otherwise, $C'$ is contained in~$C$, and the subgraphs in
  $H_1,\dots,H_{i+1}$ that are adjacent to $C'$ are a subset of the at most
  $s$ subgraphs $Q_1,\dots,Q_k,H_{i+1}$, which are pairwise adjacent since
  $F_1$ intersects all of $A_1,\dots,A_k$. In both cases, property
  $(\star)$ is maintained.

  Now assume that $k=s$. If $m\geq t$, then contracting each of
  $Q_1,\dots,Q_s,F_1,\dots,F_t$ to a single vertex gives a $K^*_{s,t}$
  minor. So we are left with the case $m\leq t-1$. Let $H_{i+1}$ be the
  subgraph of~$C$ induced by $V(F_1)\cup\dots\cup V(F_m)$. Hence $H_{i+1}$
  is induced by the union of $p_{i+1}$ paths $P_1,\dots,P_{p_{i+1}}$, where
  $p_{i+1}\leq ms\leq s(t-1)$ and each $P_j$ is a shortest path in
  $G-\Bigl(\bigl(V(H_1)\cup\dots\cup V(H_i)\bigr)\cup
  \bigl(V(P_1)\cup\dots\cup V(P_{j-1})\bigr)\Bigr)$. Consider a component
  $C'$ of $G-\bigl(V(H_1)\cup\dots\cup V(H_{i+1})\bigr)$. If $C'$ is
  disjoint from $C$, then $C'$ is a component of
  $G-\bigl(V(H_1)\cup\dots\cup V(H_i)\bigr)$ and $C'$ is not adjacent to
  $H_{i+1}$. Otherwise, $C'$ is contained in $C$, and $C'$ does not
  intersect some $A_j$ by the maximality of $m$. Thus $C'$ is adjacent to a
  subset of at most $s$ subgraphs in $Q_1,\dots,Q_s,H_{i+1}$, which are
  pairwise adjacent (since $H_{i+1}$ intersects all of $A_1,\dots,A_s$). In
  both cases, property~$(\star)$ is maintained.
\end{proof}

\subsection{Layered Treewidth and Generalised Colouring Numbers}

This subsection explores connections between layered treewidth and strong
colouring numbers. The \emph{layered width} of a tree decomposition
$(T_x:x\in V(T))$ of a graph $G$ is the minimum integer~$k$ such that, for
some layering $(V_1,\dots,V_\ell)$ of $G$, each bag $T_x$ contains at most
$k$ vertices in each layer $V_i$. The \emph{layered treewidth} of a graph
$G$ is the minimum layered width of a tree decomposition of $G$. Layered
treewidth was introduced independently by \citet{DMW17} and
\citet{Shahrokhi13}. Applications of layered treewidth include
nonrepetitive graph colouring \citep{DMW17}, queue and track layouts
\citep{DMW17}, graph drawing \citep{DMW17,BDDEW}, book embeddings
\citep{DF16}, and intersection graph theory \citep{Shahrokhi13}.

\begin{lem}\label{LayeredTreewidthSCOL}\mbox{}\\*
  Every graph $G$ with layered treewidth $k$ satisfies
  $\scol_r(G)\leq k(2r+1)$.
\end{lem}

\preproof\begin{proof}
  Let $(T_x:x\in V(T))$ be a tree decomposition of $G$ with layered width
  $k$ with respect to some layering $(V_1,\dots,V_\ell)$ of $G$. Root $T$
  at an arbitrary node $r$. For each vertex $v$ of~$G$, let~$h(v)$ be the
  node $x$ of $T$ closest to $r$ in $T$ such that $v\in T_x$, let
  $S_v$ be the subtree of~$T$ rooted at $h(v)$, and let $X_v=V(S_v)$.

  Let $\preceq^-$ be the partial ordering of $V(G)$ such that if
  $X_w\subsetneq X_v$, then $v\preceq^-w$. Let $\preceq$ be any linear
  ordering of $V(G)$ that is an extension of $\preceq^-$. By the definition
  of a tree decomposition, for any edge $vw$ with $v\preceq w$ we have
  $X_w\subseteq X_v$ and both~$v$ and~$w$ are in~$T_{h(w)}$ (since there
  must be some bag $T_u$ with $v,w\in V(T_u)$). This also means that if
  $u_1,\ldots,u_k$ is a path and $u_i$ is minimal among $u_1,\ldots,u_k$
  with respect to $\preceq$, then $X_{u_j}\subseteq X_{u_i}$ for all
  $j\in[k]$.

  Now assume that $x\in\sreach(G,\preceq,v)$ for some $v,x\in V(G)$. Thus
  $G$ contains a path $v,y_1,\dots,y_k,x$ of length at most~$r$ with
  $x\preceq v\preceq y_i$ for all $i\in[k]$. By the above observations this
  means that $x\in V(T_{h(y_k)})$ (since $x\preceq y_k$ and $y_kx$ is an
  edge), $X_v\subseteq X_x$ (since~$x$ is minimal with respect to~$\preceq$
  on the path $v,y_1,\dots,y_k,x$), and $X_{y_k}\subseteq X_v$ (since~$v$
  is minimal with respect to~$\preceq$ on the path $v,y_1,\dots,y_k$). So
  we have $X_{y_k}\subseteq X_v\subseteq X_x$, hence~$T_{h(v)}$ is on the
  path from~$T_{h(x)}$ to~$T_{h(y_k)}$ in~$T$. Since~$x$ is in both
  $T_{h(x)}$ and $T_{h(y_k)}$, $x$ must also be in every bag in the path
  from $T_{h(x)}$ to~$T_{h(y_k)}$ in~$T$, and hence $x\in V(T_{h(v)})$.
  Moreover, if $v\in V_i$, then, since $x$ has distance at most~$r$
  from~$v$ in~$G$, we have
  $x\in V_{i-r}\cup V_{i-r+1}\cup\dots\cup V_{i+r}$. Since
  $|V(T_{h(v)})\cap V_j|\leq k$, there are at most $(2r+1)k$ such
  vertices~$x$. It follows that $\scol_r(G)\leq k(2r+1)$.
\end{proof}

While $n$-vertex planar graphs may have treewidth as large as $\sqrt{n}$,
\citet{DMW17} proved that every graph with Euler genus $g$ has layered
treewidth at most $2g+3$. (More generally, \citet{DMW17} proved that a
minor-closed class of graphs has bounded layered treewidth if and only if
it excludes some apex graph as a minor.) Then \cref{LayeredTreewidthSCOL}
implies that every planar graph $G$ satisfies
\[\scol_r(G)\leq3(2r+1).\]
This result is close to the best known result, which is
$\scol_r(G)\leq5r+1$, proved by \citet{HOQRS}. More generally, by
\cref{LayeredTreewidthSCOL}, every graph $G$ with Euler genus $g$ satisfies
\[\scol_r(G)\leq (2g+3)(2r+1).\]
Again, this result is close to the best known result, which is
$\scol_r(G)\leq (4g+5)r+2g+1$, again due to \citet{HOQRS}.

\cref{LayeredTreewidthSCOL} is also interesting because it leads to linear
bounds on $\scol_r$ for non-minor-closed classes. We give three examples.

A graph is \emph{$(g,k)$-planar} if it can be drawn on a surface of Euler
genus at most $g$ with at most~$k$ crossings on each edge. Even
$(0,1)$-planar graphs can contain arbitrarily large complete graph minors
\citep{DEW16}. Nevertheless, \citet{DEW16} proved that every $(g,k)$-planar
graph has layered treewidth at most $(4g+6)(k+1)$, and this bound is tight
up to a constant factor. Then \cref{LayeredTreewidthSCOL} implies that for
every $(g,k)$-planar graph $G$,
\[\scol_r(G)\leq (4g+6)(k+1)(2r+1).\]

Map graphs are defined as follows. Start with a graph $G_0$ embedded in a
surface of Euler genus~$g$, with each face labelled a `nation' or a `lake',
such that each vertex of $G_0$ is incident with at most $d$ nations. Define
a graph $G$ whose vertices are the nations of $G_0$, where two vertices are
adjacent in $G$ if the corresponding faces in $G_0$ share a vertex. Then
$G$ is called a \emph{$(g,d)$-map graph}. A $(0,d)$-map graph is called a
\emph{(plane) $d$-map graph}; see \citep{Chen-JGT07,DFHT05,CGP02,DEW16} for
example. It is easily seen that $(g,3)$-map graphs are precisely the graphs
of Euler genus at most $g$ \citep{CGP02,DEW16}. So $(g,d)$-map graphs
provide a natural generalisation of graphs embedded in a surface. Note that
if a vertex of $G_0$ is incident with $d$ nations, then $G$ contains $K_d$,
so map graphs can have arbitrarily large complete minors. \citet{DEW16}
proved that every $(g,d)$-map graph on $n$ vertices has layered treewidth
at most $(2g+3)(2d+1)$, and this bound is tight up to a constant factor. So
\cref{LayeredTreewidthSCOL} implies that for every $(g,d)$-map graph $G$,
\[\scol_r(G)\leq (2g+3)(2d+1)(2r+1).\]

For a set $P$ of points in the plane, the \emph{unit disc} graph $G$ of $P$
has vertex set $P$, where $vw\in E(G)$ if and only if $\dist(v,w)\leq1$
(where now $\dist(v,w)$ denotes the Eulerian distance in the plane).
\citet{BDDEW} proved that every unit disc graph with maximum clique size
$k$ has layered pathwidth, and thus layered treewidth, at most $4k$. Then
\cref{LayeredTreewidthSCOL} implies that every unit disc graph with maximum
clique size $k$ satisfies
\[\scol_r(G)\leq4k(2r+1).\]

\section{Excluded Immersions}
\label{Immersions}

This section studies the defective chromatic number of graphs excluding a
fixed immersion. A graph $G$ contains a graph $H$ as an \emph{immersion}
(also called a \emph{weak immersion}) if the vertices of~$H$ can be mapped
to distinct vertices of $G$, and the edges of $H$ can be mapped to pairwise
edge-disjoint paths in $G$, such that each edge $vw$ of $H$ is mapped to a
path in $G$ whose endpoints are the images of $v$ and~$w$. The image in $G$
of each vertex in $H$ is called a \emph{branch vertex}. A graph $G$
contains a graph $H$ as a \emph{strong immersion} if $G$ contains $H$ as an
immersion such that for each edge $vw$ of $H$, no internal vertex of the
path in $G$ corresponding to $vw$ is a branch vertex.

Inspired no doubt by Hadwiger's Conjecture, \citet{LM89} and \citet{AL03}
independently conjectured that every $K_t$-immersion-free graph is properly
$(t-1)$-colourable. Often motivated by this question, structural and
colouring properties of graphs excluding a fixed immersion have recently
been widely studied \citep{GKT15,Wollan15,DW16,DDFMMS14,DK14,DMMS13,RS10}.
The best upper bound, due to \citet{GLW17}, says that every
$K_t$-immersion-free graph is properly $(3.54t+3)$-colourable.

We prove that the defective chromatic number of $K_t$-immersion-free graphs
equals~2.

\begin{thm}\label{2DefectiveImmersion}\mbox{}\\*
  Every graph not containing $K_t$ as an immersion is $2$-colourable with
  defect $(t-1)^3$.
\end{thm}

\begin{thm}\label{2DefectiveStrongImmersion}\mbox{}\\*
  For every integer $t$, there is an integer $d$ such that every graph not
  containing $K_t$ as a strong immersion is $2$-colourable with defect $d$.
\end{thm}

Notice that immersions naturally also appear in the setting of
\emph{multigraphs}, allowing multiple edges but no loops. It is obvious
that if~$G$ is a multigraph with edge multiplicity at most~$m$, then the
results of the theorems above hold with defect $m(t-1)^3$ and $md$,
respectively. On the other hand, if every edge in a multigraph has
multiplicity $m+1$, then no two adjacent vertices get the same colour in a
colouring with defect $m$. In particular, the graph obtained by replacing
the edges in the complete graph $K_{t-1}$ by $m+1$ parallel edges does not
have $K_t$ as an immersion, but is also not $(t-2)$-colourable with defect
$m$.

We leave as an open problem to determine the clustered chromatic number of
graphs excluding a (strong or weak) $K_t$ immersion. It was observed by
both \citet{HST03} and Liu and Oum~\citep{LiuOum} that the results in
\citet{ADOV03} prove that for every $k,N$, there exists a $(4k-2)$-regular
graph~$G$ such that every $k$-colouring of~$G$ has a monochromatic
component of size at least~$N$. In other words, the clustered chromatic
number of graphs with maximum degree $\Delta$ is at least
$\bigl\lfloor\frac14(\Delta+6)\bigr\rfloor$. Since every graph with maximum
degree at most $t-2$ contains no (strong or weak) $K_t$ immersion, the
clustered chromatic number of graphs excluding a (strong or weak) $K_t$
immersion is at least $\bigl\lfloor\frac14(t+4)\bigr\rfloor$.

The proof of \Cref{2DefectiveImmersion} uses the following structure
theorem from \citet{DMMS13}. The theorem is not explicitly proved
in~\cite{DMMS13}, but can be derived easily from the proof of Theorem~1 on
page~4 of that paper.

For each edge $xy$ of a tree $T$, let $T(xy)$ and $T(yx)$ be the components
of $T-xy$, where~$x$ is in~$T(xy)$ and $y$ is in~$T(yx)$. For a tree $T$
and graph $G$, a \emph{$T$-partition} of $G$ is a partition
$\bigl(T_x\subseteq V(G):x\in V(T)\bigr)$ of~$V(G)$ indexed by the nodes of
$T$. As before, each set $T_x$ is called a \emph{bag}. Note that a bag may
be empty. For each edge $xy\in E(T)$, let
$G(T,xy)=\bigcup_{z\in V(T(xy))}T_z$ and
$G(T,yx)=\bigcup_{z\in V(T(yx))}T_z$. Let $E(T,xy)$ ($=E(T,yx)$) be the set
of edges in $G$ between $G(T,xy)$ and $G(T,yx)$. The \emph{adhesion} of a
$T$-partition is the maximum, taken over all edges~$xy$ of $T$, of
$|E(T,xy)|$. For each node $x$ of $T$, the \emph{torso of $x$} (with
respect to a $T$-partition) is the graph obtained from $G$ by identifying
$G(T,yx)$ into a single vertex for each edge $xy$ incident to $x$, deleting
resulting parallel edges and loops.

\begin{thm}[following \citet{DMMS13}]\label{WeakImmStructure}\mbox{}\\*
  For every graph $H$ with $t$ vertices and every graph $G$ that does not
  contain $H$ as an immersion, there is a tree $T$ and a $T$-partition of
  $\,G$ with adhesion less than $(t-1)^2$, such that each bag has at most
  $t-1$ vertices.
\end{thm}

A structural result similar to this theorem was proved by \citet{Wollan15}.
We also need the following lemma.





\begin{lem}\label{Bijection}\mbox{}\\*
  Let $G$ be a graph such that for some tree $T$ with vertex set $V(G)$,
  for each edge $xy$ of $\,T$, the number of edges of $\,G$ between
  $V(T(xy))$ and $V(T(yx))$ is at most $k$. Then $G$ is $2$-colourable with
  defect $k$.
\end{lem}

\preproof\begin{proof}
  We use induction on $|V(G)|$, noting that there is nothing to prove if
  $|V(G)|\leq2$. So assume $|V(G)|\geq3$. Call a vertex~$v$ of~$G$
  \emph{large} if $\deg_G(v)\geq k+1$; otherwise~$v$ is \emph{small}.

  If $G$ has no large vertices, then every 2-colouring of $G$ has defect
  $k$. Now assume that $G$ has some large vertex. Thus there is an edge
  $uv$ of $T$ such that $u$ is large and $u$ is the only large vertex in
  $V(T(uv))$. Set $a=|V(T(uv))|$. Suppose that every vertex in
  $V(T(uv))\setminus\{u\}$ has a neighbour in~$G$ in $V(T(vu))$. Since~$u$
  has at least $k+1-(a-1)$ neighbours outside~$V(T(uv))$, the number of
  edges between $V(T(uv))$ and $V(T(vu))$ is at least
  $\bigl(k+1-(a-1)\bigr)+(a-1)=k+1$, a contradiction.

  So there is a vertex $w\in V(T(uv))\setminus\{u\}$ with $N_G(w)\subseteq
  V(T(uv))$. Note that~$w$ is small. Let~$wz$ be an edge in $T$. Form the
  graphs $G'$ and $T'$ respectively from $G$ and $T$ by identifying~$w$ and
  $z$ (deleting loops and parallel edges). For each edge $xy$ of~$T'$, the
  number of edges of~$G'$ between $V(T'(xy))$ and $V(T'(yx))$ is still at
  most~$k$. Hence by induction, $G'$ has a $2$-colouring with defect~$k$.
  This colouring gives a $2$-colouring with defect~$k$ of all vertices
  of~$G$ except~$w$. Since all vertices in $V(T(uv))$ except~$u$ are small,
  $u$ is the only possible large neighbour of $w$. Give~$w$ the colour
  different from $u$. As all other neighbours of~$w$ are small, the
  monochromatic degree can increase only for small vertices. Thus the
  defect is at most $k$, as required.
\end{proof}

Now we are ready to prove our $2$-colour result for graphs excluding an
immersion.

\begin{proof}[Proof of \Cref{2DefectiveImmersion}.]
  By \cref{WeakImmStructure}, there is a tree $T$ and a $T$-partition of
  $G$ with adhesion at most $(t-1)^2-1$, such that each bag has at most
  $t-1$ vertices. Let $Q$ be the graph with vertex set $V(T)$, where $xy\in
  E(Q)$ whenever there is an edge of $G$ between $T_x$ and $T_y$. Any one
  edge of~$Q$ corresponds to at most $t-1$ edges in $G$. By
  \cref{Bijection}, the graph $Q$ is $2$-colourable with defect
  $(t-1)^2-1$. Assign to each vertex $v$ in $G$ the colour assigned to the
  vertex $x$ in~$Q$ with $v\in T_x$. Since at most $t-1$ vertices of $G$
  are in each bag, $G$ is $2$-coloured with defect at most
  $(t-1)\cdot\bigl((t-1)^2-1\bigr)+(t-2)<(t-1)^3$.
\end{proof}

To prove our result for strong immersions, we employ the following more
involved structure theorem of \citet{DW16}.

\begin{thm}[\citet{DW16}]\label{StrongImmersionStructure}\mbox{}\\*
  For every integer $t$, there is an integer $\alpha$ such that for every
  graph $G$ that does not contain~$K_t$ as a strong immersion, there is a
  tree $T$ and a $T$-partition of $\,G$ with adhesion at most $\alpha$ such
  that the following holds. For each node~$x$ of~$\,T$ with torso $S_x$, if
  $W_x$ is the set of vertices in~$S_x$ with degree at least~$\alpha$, then
  there is a subset $A_x\subseteq W_x$ of size at most $\alpha$ such that
  $W_x\setminus A_x$ can be enumerated $\{x_1,\dots,x_p\}$ and $V(S_x-W_x)$
  can be partitioned $B_0,B_1,\dots,B_p$ (allowing $B_j=\varnothing$), such
  that:
  \vspace*{-2.25ex}
  \begin{enumerate}[(1)]
    \setlength\itemsep{0ex}
  \item each vertex $v\in A_x$ is adjacent to at most $\alpha$ of
    $\,B_0,B_1,\dots,B_p$ and adjacent to at most $\alpha$ vertices in
    $W_x\setminus A_x$; and
  \item for each $i\in[p]$, there are at most $\alpha$ edges
    between $B_0\cup\dots\cup B_{i-1}\cup\{x_1,\dots,x_{i-1}\}$ and
    $B_i\cup\dots\cup B_p\cup\{x_{i+1},\dots,x_p\}$.
  \end{enumerate}
\end{thm}

We actually only need the following corollary of
\cref{StrongImmersionStructure}.

\begin{cor}\label{StrongImmersionStructureCorollary}\mbox{}\\*
  For every integer $t$, there is an integer $\alpha$ such that for every
  graph $G$ that does not contain~$K_t$ as a strong immersion, there is a
  tree $T$ and $T$-partition of $\,G$ with adhesion at most $\alpha^2$ such
  that for each node~$x$ of~$\,T$ with torso $S_x$, $S_x[T_x]$ has degree
  at most $3\alpha+2$.
\end{cor}

\preproof\begin{proof}
  Consider a tree~$T$ and $T$-partition of $G$ in accordance with
  \cref{StrongImmersionStructure}. Consider a node $x$ of $T$ with torso
  $S_x$. We use the notation from the theorem.

  Consider a vertex $v\in W_x$. If $v\in A_x$, then~$v$ has at most
  $|A_x|-1<\alpha$ neighbours in $A_x$ and at most $\alpha$ neighbours in
  $W_x\setminus A_x$, and thus has less than $2\alpha$ neighbours in $W_x$.
  If $v\in W_x\setminus A_x$, then $v=x_i$ for some $i\in[p]$. Then~$v$ has
  at most $|A_x|\leq\alpha$ neighbours in $A_x$. Furthermore, there are at
  most $\alpha$ edges between $\{x_1,\dots,x_{i-2}\}$ and
  $\{x_i,\dots,x_p\}$, at most $\alpha$ edges between $\{x_1,\dots,x_i\}$
  and $\{x_{i+2},\dots,x_p\}$, and at most~$2$ edges between $x_i$ and
  $\{x_{i-1},x_{i+1}\}$. Thus $v$ has at most $3\alpha+2$ neighbours in
  $W_x$. Hence $S_x[W_x]$ has maximum degree at most $3\alpha+2$.

  Apply the following operation for each vertex $v\in T_x\setminus W_x$,
  for each node $x$ of $T$ with $|T_x|\geq2$. Since $v\not\in W_x$, the
  degree of $v$ in $S_x$ is at most $\alpha-1$. Since there are at most
  $\alpha$ edges from~$G$ between $T_x$ and each contracted vertex in
  $S_x$, $v$ has degree at most $(\alpha-1)\alpha<\alpha^2$ in $G$. Now
  delete $v$ from $T_x$, add a new node $y$ in $T$ adjacent only to $x$,
  and define $T_y=\{v\}$. Note that the number of edges between $T_y$ and
  $G-T_y$ is less than $\alpha^2$, and the torso of $y$ is isomorphic to
  $K_2$ (hence has degree one). Finally, the torso of~$x$ hasn't changed,
  since the contraction of the single-vertex node $T_y$ just gives the
  vertex~$v$ again. In particular, the degree of~$v$ in the torso of~$x$ is
  still at most $\alpha-1$, and hence~$W_x$ also hasn't changed.

  After having applied the operation from the previous paragraph as long a
  possible, we obtain a tree-partition of~$G$ on a tree $T'$ with adhesion
  at most $\alpha^2$. Moreover, for each node $x$ of~$T'$ we have
  $|T'_x|=1$, and then $S'_x[T'_x]$ has degree zero, or
  $T'_x\subseteq W'_x$ and $S'_x[W'_x]$ has degree at most $3\alpha+2$. We
  immediately get that $S'_x[T'_x]$ has degree at most $3\alpha+2$ as well.
\end{proof}

Now we are ready to prove our $2$-colouring result for graphs excluding a
strong immersion.

\begin{proof}[Proof of \Cref{2DefectiveStrongImmersion}.]
  By \cref{StrongImmersionStructureCorollary}, there is an
  integer~$\alpha$, a tree $T$ and a $T$-partition of $G$ with adhesion at
  most $\alpha^2$, such that for each node $x$ of~$T$ with torso $S_x$,
  $S_x[T_x]$ has degree at most $3\alpha+2$. Let $Q$ be the graph with
  vertex set $V(T)$, where $xy\in E(Q)$ whenever there is an edge of $G$
  between $T_x$ and~$T_y$. By \cref{Bijection} with $k=\alpha^2$, the graph
  $Q$ is $2$-colourable with defect~$\alpha^2$. Assign to each vertex~$v$
  in $G$ the colour assigned to the vertex $x$ in $Q$ with $v\in T_x$.

  If $v\in V(T_x)$, then every edge~$vw$ in~$G$ with $w\not\in T_x$ gives
  rise to an edge in~$Q$. Since the adhesion is at most~$\alpha^2$, any one
  edge of~$Q$ corresponds to at most~$\alpha^2$ edges in~$G$. As the
  monochromatic degree of~$x$ in~$Q$ is at most~$\alpha^2$, this means
  that~$v$ has at most~$\alpha^4$ neighbours outside~$T_x$ with the same
  colour. Adding the at most $3\alpha+2$ neighbours of~$v$ in~$T_x$, we
  obtain that the monochromatic degree of~$v$ in~$G$ is at most
  $\alpha^4+3\alpha+2$.
\end{proof}

\subsection*{Notes and acknowledgements}

\vspace{-1.5ex}
After publication of the first version of this paper, \citet{DN17} proved
that every graph embeddable in a fixed surface is $4$-colourable with
bounded clustering (cf.\ the comments after \cref{K3tColouring}). They also
gave an alternative proof that the clustered chromatic number of
$K_t$-minor-free graphs is at most $2t-2$ (cf.\ \cref{MonoComponents}), and
announced that in a sequel they will prove that the clustered chromatic
number of $K_t$-minor-free graphs equals $t-1$.

Thanks to Jacob Fox who first observed that bounded degree graphs give
lower bounds on the clustered chromatic number of graphs excluding a fixed
immersion.

The authors also like to thank an anonymous referee for pointing out some errors in earlier versions of this paper.

  \let\oldthebibliography=\thebibliography
  \let\endoldthebibliography=\endthebibliography
  \renewenvironment{thebibliography}[1]{%
    \begin{oldthebibliography}{#1}%
      \setlength{\parskip}{0.25ex}%
      \setlength{\itemsep}{0.25ex}%
  }%
  {%
    \end{oldthebibliography}%
  }

\def\soft#1{\leavevmode\setbox0=\hbox{h}\dimen7=\ht0\advance \dimen7
  by-1ex\relax\if t#1\relax\rlap{\raise.6\dimen7
  \hbox{\kern.3ex\char'47}}#1\relax\else\if T#1\relax
  \rlap{\raise.5\dimen7\hbox{\kern1.3ex\char'47}}#1\relax \else\if
  d#1\relax\rlap{\raise.5\dimen7\hbox{\kern.9ex \char'47}}#1\relax\else\if
  D#1\relax\rlap{\raise.5\dimen7 \hbox{\kern1.4ex\char'47}}#1\relax\else\if
  l#1\relax \rlap{\raise.5\dimen7\hbox{\kern.4ex\char'47}}#1\relax \else\if
  L#1\relax\rlap{\raise.5\dimen7\hbox{\kern.7ex
  \char'47}}#1\relax\else\message{accent \string\soft \space #1 not
  defined!}#1\relax\fi\fi\fi\fi\fi\fi}

\end{document}